\documentclass[12pt,reqno,a4paper]{amsart}

\usepackage{a4}
\usepackage{amsmath, amssymb}
\usepackage{amsthm}

\usepackage{graphicx}                      
\newcommand{\color}[2][{}]{}        

\usepackage{bbm}         
\usepackage{mathrsfs} 
\renewcommand\mathcal\mathscr  

\usepackage{accents}   

\usepackage{slashed}    

\usepackage{stmaryrd}    

\usepackage{diagrams}               

\theoremstyle{plain}            
\newtheorem{theorem}{Theorem}[section]

\newtheorem{lemma}[theorem]{Lemma}

\theoremstyle{definition}       
\newtheorem{definition}[theorem]{Definition}

\newtheorem{example}[theorem]{Example}

\theoremstyle{remark}           
\newtheorem{remark}[theorem]{Remark}

\newcommand{\Sec}[1]{Section~\ref{sec:#1}}

\newcommand{\ExS}[2]{Examples~\ref{ex:#1}--\ref{ex:#2}}
\newcommand{\Exenum}[2]{Example~\ref{ex:#1}~(\ref{#2})}
\newcommand{\Exenums}[3]{Example~\ref{ex:#1}~(\ref{#2})--(\ref{#3})}
\newcommand{\Eq}[1]{Eq.~\eqref{eq:#1}}
\newcommand{\Eqs}[2]{Eqs.~\eqref{eq:#1} and~\eqref{eq:#2}}

\newcommand{\Thm}[1]{Theorem~\ref{thm:#1}}
\newcommand{\Thmenum}[2]{Theorem~\ref{thm:#1}~(\ref{#2})}

\newcommand{\Lem}[1]{Lemma~\ref{lem:#1}}
\newcommand{\Lemenum}[2]{Lemma~\ref{lem:#1}~(\ref{#2})}

\newcommand{\Rem}[1]{Remark~\ref{rem:#1}}
\newcommand{\Remenum}[2]{Remark~\ref{rem:#1}~(\ref{#2})}
\newcommand{\Def}[1]{Definition~\ref{def:#1}}

\newcommand{\Defenum}[2]{Definition~\ref{def:#1}~(\ref{#2})}
\newcommand{\Defenums}[3]{Definition~\ref{def:#1}~(\ref{#2})--(\ref{#3})}

\numberwithin{equation}{section}

\DeclareMathOperator{\dist}   {dist}
\DeclareMathOperator{\dom}    {dom}
\DeclareMathOperator{\ran}    {ran}

\DeclareMathOperator{\vol}    {vol}

\DeclareMathOperator{\const}  {const}

\newcommand{\spec}[2][{}]   {\sigma_{\mathrm{#1}}(#2)}

\def\XXint#1#2#3{{\setbox0=\hbox{$#1{#2#3}{\int}$}
     \vcenter{\hbox{$#2#3$}}\kern-.5\wd0}}

\newlength{\maxbreite}%
\newlength{\maxhoehe}%
\newlength{\maxtiefe}%

\newcommand{\stelldrueber}[3][0pt]{
  \settowidth{\maxbreite}{#3}%
  \settoheight{\maxhoehe}{#3}%
  \settodepth{\maxtiefe}{#2}%
  \addtolength{\maxhoehe}{\maxtiefe}%
  {\makebox[\maxbreite]{\raisebox{\maxhoehe}{\hspace{#1}#2}}%
  \makebox[0pt][r]{#3}}%
}

\newcommand{\overcirc}[1]       
{\stelldrueber[.45ex]{$\scriptscriptstyle\circ$}{${#1}$}}

\newcommand{\R}{\mathbb{R}} 
\newcommand{\C}{\mathbb{C}} 
\newcommand{\Z}{\mathbb{Z}} 

\newcommand{\eps}{\varepsilon} 
\renewcommand{\phi}{\varphi}   
\newcommand{\e}{\mathrm e}  
\newcommand{\im}{\mathrm i} 
\DeclareMathOperator{\dd}    {d\!} 

\newcommand{\wt}{\widetilde}           
\newcommand {\qf}[1]{\mathfrak{#1}}    

\newcommand{\HS}{\mathcal H}           

\newcommand{\Sobsymb} {\mathsf H}      
\newcommand{\Contsymb} {\mathsf C}     
\newcommand{\Lsymb}    {\mathsf L}     
\newcommand{\lsymb}    {\ell}          
\newcommand{\Lsqrsymb}    {\mathsf L_2}     
\newcommand{\lsqrsymb}    {\ell_2}          

\newcommand{\Cont}[2][{}]{\Contsymb^{#1}({#2})}

\newcommand{\Lsqr}[2][{}]{\Lsymb_2^{#1}({#2})} 
 
\newcommand{\lsqr}[2][{}]{\lsymb_2^{#1}({#2})}   

\newcommand{\Sob}[2][1]{\Sobsymb^{#1}({#2})} 
\newcommand{\Sobn}[2][1]{\Sobsymb_\circ^{#1}({#2})}
 
\newcommand{\norm}[2][{}]{\|{#2}\|_{{#1}}}    
\newcommand{\normsqr}[2][{}]{\|{#2}\|^2_{#1}} 

\newcommand{\iprod}[3][{}]{\langle{#2},{#3}\rangle_{#1}}  
\newcommand{\bigiprod}[3][{}]{\bigl\langle{#2},{#3}\bigr\rangle_{#1}}
\newcommand{\Bigiprod}[3][{}]{\Bigl\langle{#2},{#3}\Bigr\rangle_{#1}}

\newcommand{\set}[2]{\{ \, #1 \, | \, #2 \, \} } 
\newcommand{\bigset}[2]{\bigl\{ \, #1 \, \bigl|\bigr. \, #2 \, \bigr\} }
\newcommand{\Bigset}[2]{\Bigl\{ \, #1 \, \Bigl|\Bigr. \, #2 \, \Bigr\} }

\newcommand{\map}[3]{ #1 \colon #2 \longrightarrow #3 } 

\newcommand{\bd}  {\partial}                
\newcommand{\clo}[1]{\overline{{#1}}} 

\newcommand{\dcup}{\mathbin{\mathaccent\cdot\cup}}

\DeclareMathOperator*{\bigdcup}{\mathaccent\cdot{\bigcup}}

\newcommand{\restr}[1]{{\restriction}_{#1}} 

\newcommand{\conj}[1]{\overline {{#1}}}       
\newcommand{\orth}{\bot}                    

\newcommand{\1}{\mathbbm 1}                    
\newcommand{\und}{\quad\text{and}\quad}

\newcommand{\orient}[1]{\accentset{\curvearrowright}{#1}} 

\newcommand{\Neu}{{\mathrm N}}              
\newcommand{\Dir}{{\mathrm D}}              
\newcommand{\laplacian}[2][{}]{\Delta_{{#2}}^{{#1}}} 
\newcommand{\laplacianD}[1]{\laplacian[\Dir]{#1}} 
\newcommand{\laplacianN}[1]{\laplacian[\Neu]{#1}} 

\newcommand{\vxeps}{{\eps,v}}

\usepackage{amsbsy}  

\newcommand{\dbar}[1]{\leavevmode\raise0.6ex\hbox{--}\kern-0.7em #1}

\newcommand{\mc}{\mathcal}
\newcommand{\ul}{\underline}
\newcommand{\orul}[1]{\orient {\underline{#1}}}

\newcommand{\mbE}{\mathbb E} 

\newcommand{\Graph} X  
\newcommand{\Gmax}{\mc G^{\max}}

\newcommand{\Forms}{\Lambda}

\newcommand{\oplusmerge}{\stackrel{\curlywedgeuparrow}{\oplus}}
\newcommand{\oplussplit}{\stackrel{\curlyveeuparrow}{\oplus}}

\newcommand{\de}   {\mathord {\mathrm d}}         
\DeclareMathOperator{\ind}   {ind}      

\newcommand{\dde}{\mathsf{d}}
\declareslashed{}{\text{-}}{-0.13}{0.4}{D}
\newcommand{\De}{\mathsf{D}}

\newcommand{\dlaplacian}[2][{}]{\pmb{\triangle}_{{#2}}^{{#1}}}

\newcommand{\wtdlaplacian}[2][{}]{{%
  \wt \triangle_{{#2}}^{{#1}}%
}}                              

\newcommand{\stand}{\mathrm{std}}   
\newcommand{\magn}{\mathrm{mag}}

\newcommand{\bigoplussqr}{\bigoplus}
\newcommand{\orientul}[1]{\underline {\orient {#1}}} 
\newcommand{\orientmap}{(\hspace{0.2ex}\orient \cdot \hspace{0.1ex})}

\newcommand{\Sobx}[3][1]{\Sobsymb_{{#2}}^{#1}({#3})} 

\newcommand{\embmap}[3]{ #1 \colon #2 \hookrightarrow #3 } 

\begin{document}
\title[First order approach and index theorems for graphs]{First order
  approach and index theorems for discrete and metric graphs}

\author{Olaf Post}      
\address{Institut f\"ur Mathematik,
         Humboldt-Universit\"at zu Berlin,
         Rudower Chaussee~25,
         12489 Berlin,
         Germany}
\email{post@math.hu-berlin.de}
\date{\today}




\begin{abstract}
  The aim of the present paper is to introduce the notion of first
  order (supersymmetric) Dirac operators on discrete and metric
  (``quantum'') graphs.  In order to cover all self-adjoint boundary
  conditions for the associated metric graph Laplacian, we develop
  systematically a new type of discrete graph operators acting on a
  decorated graph. The decoration at each vertex of degree~$d$ is
  given by a subspace of $\C^d$, generalising the fact that a function
  on the standard vertex space has only a scalar value.

  We develop the notion of exterior derivative, differential forms,
  Dirac and Laplace operators in the discrete and metric case, using a
  supersymmetric framework. We calculate the (supersymmetric) index of
  the discrete Dirac operator generalising the standard index formula
  involving the Euler characteristic of a graph. Finally, we show that
  the corresponding index for the metric Dirac operator agrees with
  the discrete one.
\end{abstract}

\maketitle

%
\section{Introduction}
\label{sec:intro}
%

In the last years, many attention has been payed in the analysis of
metric graph Laplacians, i.e., operators acting as second order
differential operators on each edge considered as one-dimensional
space, with suitable (vertex) boundary conditions turning the
Laplacian into a self-adjoint (unbounded) operator. In most of the
works, the \emph{second order} operator is the starting object for the
analysis. For more details on Laplacians on metric graphs, also
labelled as ``quantum graphs'', we refer to the
articles~\cite{kostrykin-schrader:06,kuchment:04,kuchment:05} and the
references therein.

In this paper whereas, we want to introduce the metric graph
Laplacians with general (non-negative) vertex boundary conditions via
\emph{first order} operators, namely via an exterior derivative
analogue as in differential geometry. As a by-product, we obtain a new
type of discrete graph operators acting on a decorated graph. The
decoration at each vertex~$v$ of degree~$\deg v$ is given by a
subspace of~$\C^{\deg v}$, generalising the fact that a function~$F
\in \lsqr V$ on the standard vertex space on~$V$ has only a scalar
value~$F(v) \in \C$. In addition, we introduce the notion of a
discrete exterior derivative, a discrete Dirac and Laplace operator
and show an index theorem generalising the standard index formula
involving the Euler characteristic of a graph (cf.\
\Thm{index.discr}).

In a second part, we define exterior derivatives, Dirac and Laplace
operators on a (continuous) metric graph and relate their kernels with
the appropriate discrete objects and show that the index agrees with
the index of the discrete setting (cf.\ \Thm{index}). 

We introduce all Laplacians in a supersymmetric setting, i.e., by
appropriate ``exterior derivatives'' mimicking the corresponding
notion for manifolds. The advantage is the simple structure of these
operators; and the use of the abstract supersymmetric setting, e.g.,
the spectral equality of the Laplacian defined on even and odd
``differential forms'' (cf.~\Lem{susy}).

Index formulas may be used in order to decide whether a metric graph
$\Graph_0$ with Laplacian $\laplacian{\Graph_0}$ occurs as limit of a
``smooth'' space, i.e., a manifold or an open neighbourhood $X_\eps$
of $\Graph_0$ together with a natural Laplacian
$\laplacian{X_\eps}$. If $X_\eps$ is homotopy-equivalent to $\Graph_0$
then their Euler characteristics agree, and correspondingly,
appropriately defined indices for the operators on $X_\eps$ and
$\Graph_0$ must agree if the operators converge. We comment on this
observation in \Sec{mg.smooth}.

Spectral graph theory is an active area of research.  We do not
attempt to give a complete overview here. Results on spectral theory
of discrete or combinatorial Laplacians can be found e.g.\
in~\cite{dodziuk:84,mohar-woess:89,colin:98,chung:97}. For continuous
(quantum) graph Laplacians we mention the
works~\cite{roth:84,nicaise:87,kostrykin-schrader:99,%
  harmer:00b,kostrykin-schrader:03,kuchment:04,%
  friedman-tillich:pre04, kuchment:05,%
  kostrykin-schrader:06,pankrashkin:06a,hislop-post:pre06}.  In
particular, a heat equation approach for the index formula for certain
metric graph Laplacian (with energy-independent scattering matrix) can
be found in~\cite{kps:pre07}. In particular, when submitting this
work, we learned about a related work on index formulas on quantum
graphs proven in a direct way (not using our discrete exterior
calculus) by Fulling, Kuchment and Wilson~\cite{fkw:pre07}.
Prof.~Fulling announced the results in a talk at the Isaac Newton
Institute (INI) in Cambridge~\cite{fulling.talk:07} where also the
first order factorisation of the standard quantum graph Laplacian
appears.

When submitting this work,  the work~\cite{fkw:pre07} where a similar index formula
for quantum graphs is proven in a direct way (not using our discrete
exterior calculus).

The paper is organised as follows: In the next subsection, we start
with a motivating example of standard boundary conditions in order to
illustrate the basic results and ideas. In \Sec{susy}, we develop the
abstract setting of supersymmetry. In \Sec{vx.sp}, we define a
generalisation for the discrete vertex space $\lsqr V$, namely,
general vertex spaces. In \Sec{op.vx.sp} we generalise the notion of
the coboundary operator (``exterior derivative''), Dirac and Laplace
operators in this context. In \Sec{ind.discr} we calculate the index
of the discrete Dirac operator for general vertex spaces and
generalise the below discrete Gau{\ss}-Bonnet
formula~\eqref{eq:gauss.bonnet}. In \Sec{quantum} we develop the
theory of ``exterior derivatives'' on a metric graph and introduce the
corresponding notion of Dirac and Laplace operators. In particular, we
cover all self-adjoint boundary conditions leading to a non-negative
Laplacian. Finally, in \Sec{kernel} we show that the discrete and
continuous Laplacians agree at the bottom of the spectrum, i.e., the
index formula~\eqref{eq:kernel.std} for the general case. We conclude
with a series of examples showing how an index formula can be used to
find ``smooth'' approximations of metric graph Laplacians.

\subsection*{Acknowledgements}
The author would like to thank the organisers of the programme
``Analysis on graphs and its applications'' at the Isaac Newton
Institute (INI) in Cambridge for the kind invitation. The very
inspiring atmosphere and many discussions led to this work.

\subsection{The standard case}
\label{sec:std}
In order to motivate our abstract setting, we start with the standard
Laplacian in the discrete and continuous setting. Details can be found
in the subsequent sections.  Let $\Graph=(V,E)$ be an oriented graph
with $V$ the set of vertices and $E$ the set of edges $e$, where we
denote the initial vertex by $\bd_-e$ and the terminal vertex by
$\bd_+$. Denote by $\lsqr V$ the standard vertex space with weight
$\deg v$, the degree of the vertex $v$. We consider a (scalar)
function in $\lsqr V$ as a ``$0$-form''. The \emph{coboundary
  operator} or \emph{(discrete) exterior derivative} is defined as
\begin{equation*}
  \map {\dde} {\lsqr V}{\lsqr E}, \qquad
  (\dde F)_e = F(\bd_+ e) - F(\bd_- e)
\end{equation*}
mapping $0$-forms into $1$-forms with adjoint operator
\begin{equation*}
  \map {\dde^*} {\lsqr E} {\lsqr V}, \qquad
  (\dde \eta)(v) = \frac 1 {\deg v}
            \sum_{e \in E_v} \orient \1_e(v) \eta_e
\end{equation*}
where $\1_e(v)=\pm 1$ if $v=\bd_\pm e$ and $E_v$ is the set of edges 
adjacent to $v$. We call the operator
\begin{equation}
  \label{eq:def.dirac.std}
  \De (F \oplus \eta) = \dde^* \eta \oplus \dde F,
    \qquad \text{i.e.,} \qquad \De \cong 
  \begin{pmatrix}
    0 & \dde^* \\ \dde & 0
  \end{pmatrix}
\end{equation}
the associated \emph{Dirac} operator on $\lsqr {\Forms \Graph} :=
\lsqr V \oplus \lsqr E$. The associated Laplacian is defined as
$\dlaplacian {\Forms \Graph} := \De^2$, and in particular, its
component on $0$-forms, i.e., on $\lsqr V$ is the standard Laplacian
of discrete graph theory, namely
\begin{equation}
  \label{eq:std.lap0}
  (\dde^* \dde F)(v) = (\dlaplacian[0] \Graph F)(v) 
  = \frac 1 {\deg v} \sum_{e \in E_v} 
    \bigl( F(v) - F(v_e) \bigr) 
\end{equation}
where $v_e$ denotes the vertex opposite to $v$ on $e \in E_v$.  For a
\emph{finite} graph $\Graph$, we define the \emph{index} of $\De$ as
\begin{equation}
  \label{eq:def.ind.std}
  \ind \De := \dim \ker \dde - \dim \ker \dde^*,
\end{equation}
i.e., the index of $\De$ is the Fredholm index of $\dde$.  It is a
classical result from cohomology theory, that the Fredholm-index of
the coboundary operator $\dde$ equals the \emph{Euler characteristic}
$\chi(\Graph) := |V|-|E|$, namely,
\begin{equation}
\label{eq:ind.std}
  \ind \De = \chi(\Graph).
\end{equation}
If we define the \emph{curvature} at the vertex~$v \in V$ as
  \begin{equation}
    \label{eq:curv}
    \kappa(v) := 1 - \frac 12 \deg v,
  \end{equation}
we can interprete the formula~\eqref{eq:ind.std} as a ``discrete
Gau{\ss}-Bonnet'' theorem, namely
\begin{equation}
  \label{eq:gauss.bonnet}
  \ind \De = \sum_{v \in V} \kappa(v)
\end{equation}
using the classical formula $2|E| = \sum_{v \in V} \deg v$.  Note
that~$\kappa(v) < 0$ iff~$\deg v \ge 3$. 

Considering $X$ as a metric graph, our basic Hilbert space is $\Lsqr
\Graph$ (cf.~\eqref{eq:lsqr.qg}). On the metric graph, 
we consider the ``exterior'' derivative 
\begin{equation*}
  \map {\de} {\dom \de} {\Lsqr \Graph}, \qquad \de f = f'=\{f_e'\}_e
\end{equation*}
where $\dom \de = \Sobx \max \Graph \cap \Cont \Graph$ is the Sobolev
space of functions \emph{continuous} at each vertex. Its
$\Lsqrsymb$-adjoint is
\begin{equation*}
  \map {\de} {\dom \de^*} {\Lsqr \Graph}, \qquad \de g = -g'= \{-g_e'\}_e
\end{equation*}
with $g \in \dom \de^*$ iff
\begin{equation}
  \label{eq:kirchhoff}
  \sum_{e \in E_v} \orient g_e(v) = 0,
\end{equation}
where $\orient g_e(v)$ is the \emph{oriented} evaluation at $v$ (see
\Eq{sign}). As before, we can define a Dirac operator $D$ on $\Lsqr
\Graph \oplus \Lsqr \Graph$ and the
associated Laplacian $\laplacian {\Forms \Graph}$ such that its
$0$-form component is
\begin{equation*}
  \laplacian [0] \Graph f := \de^* \de f = -f'' = \{-f''_e\}_e
\end{equation*}
with domain
\begin{equation*}
  \dom \laplacian [0] \Graph =
  \bigset {f \in \dom \de}  {f' \in \dom \de^*},
\end{equation*}
i.e., the \emph{standard} Laplacian on a metric graph with functions
continuous at each vertex and the Kirchoff sum condition for the
derivative at each vertex. Although the $0$- and $1$-forms are
formally the same, they differ in their interpretation: We consider
$0$-forms as \emph{scalar} functions, whereas a $1$-form is a
\emph{vector-field} with orientation. Then the Kirchhoff sum
condition~\Eq{kirchhoff} is just a ``flux'' conservation for the flux
generated by the ``vector field'' $f'$.

Again, we define the index $\ind \De$ of the metric graph Dirac
operator $\De$ as the Fredholm-index of $\de$, i.e. in the same way as
in \Eq{def.ind.std} and one of our main results in this setting
(cf.~\Thm{index}) is
\begin{equation}
  \label{eq:kernel.std}
  \ker D \cong \ker \De \und
  \ind D = \ind \De (=|V|-|E|),
\end{equation}
i.e., an isomorphism between the kernels of the discrete and
continuous case.

We want to generalise the above setting to quantum graph Laplacians
with \emph{general} self-adjoint operators $\laplacian \Graph$ (such
that $\laplacian \Graph \ge 0$) and derive a similar index formula.

%
\subsection{Supersymmetry}
\label{sec:susy}
%
Before defining several operators on a graph, we collect common
features shared by several operators. Since in our cases we only
define~$p$-forms for~$p \in \{0,1\}$, we can identify forms of even
and odd degree with the cases~$p=0$ and~$p=1$, respectively.
\begin{definition}
  \label{def:susy}
  Let~$\HS = \HS_0 \oplus \HS_1$ be a Hilbert space and~$\map \de
  {\dom \de} {\HS_1}$ a closed operator with~$\dom \de \subset \HS_0$
  ($\de$ may be bounded, in this case we have~$\dom \de =\HS_0$).
  Then we say that~$\de$ has \emph{supersymmetry} or that $\de$ is an
  \emph{exterior derivative}. A \emph{$p$-form} is an element
  in~$\HS_p$. Furthermore, we define the associated \emph{Dirac}
  operator as
  \begin{equation*}
    D(f_0 \oplus f_1) = \de^* f_1 \oplus \de f_0,
       \quad \text{i.e.,} \quad
    D \cong
    \begin{pmatrix}
      0 & \de^* \\ \de & 0
    \end{pmatrix}
  \end{equation*}
  with respect to the decomposition~$\HS=\HS_0 \oplus \HS_1$. The
  associated \emph{Laplacian} is given by~$\Delta := D^2$. In
  particular,
  \begin{equation*}
    \Delta \cong
    \begin{pmatrix}
      \Delta_0 & 0 \\ 0 & \Delta_1
    \end{pmatrix},
  \end{equation*}
  where~$\Delta_0 = \de^* \de$ and~$\Delta_1 = \de \de^*$ on their
  natural domains.
\end{definition}
Clearly,~$\Delta$ and~$\Delta^p$ are closed, non-negative operators.
Note that $\ker \de = \ker \Delta_0$ and $\ker \de^* = \ker \Delta_1$.

We denote the spectral projection of~$\Delta_p$ by~$\1_B(\Delta_p)$.
We have the following results on the spectrum away from~$0$:
\begin{lemma}
  \label{lem:susy}
  Assume that~$\de$ has supersymmetry and that~$B \subset [0,\infty)$
  is a bounded Borel set. Then
  \begin{equation*}
    \de   \1_B(\Delta_0) = \1_B(\Delta_1) \de \und
    \de^* \1_B(\Delta_1) = \1_B(\Delta_0) \de^*.
  \end{equation*}
  Furthermore, if~$0$ is \emph{not} contained in~$B$, then
  \begin{equation*}
    \map {\de}   {\1_B(\Delta_0)(\HS_0)} {\1_B(\Delta_1)(\HS_1)} \und
    \map {\de^*} {\1_B(\Delta_1)(\HS_1)} {\1_B(\Delta_0)(\HS_0)}
  \end{equation*}
  are isomorphisms. In particular,
  \begin{equation*}
    \dim \1_B(\Delta_0) = \dim \1_B(\Delta_1) \und
    \spec {\Delta_0} \setminus \{0\} =
    \spec {\Delta_1} \setminus \{0\},
  \end{equation*}
  i.e., the spectra of~$\Delta_0$ and~$\Delta_1$ away from~$0$ agree
  including multiplicity.
\end{lemma}
\begin{proof}
  The first assertion follows from~$\de \phi(\de^*\de) = \phi(\de
  \de^*) \de$, first for polynomials~$\phi$, then for
  functions~$\phi(\lambda)=(\lambda+1)^{-k}$,~$k \ge 1$, and finally
  by the spectral calculus also for (fast enough decaying) continuous
  and measurable functions. The second assertion follows since~$\ker
  \de = \ker \Delta_0=0$ and~$\ker \de^* = \ker \Delta_1=0$. The last
  statement is a simple consequence of the isomorphisms.
\end{proof}

We have the following result, an abstract version of the Hodge
decomposition:
\begin{lemma}
  \label{lem:hodge}
  Assume that $\de$ has supersymmetry and that the associated Dirac
  operator $D$ has a spectral gap at $0$, i.e., $\dist(0,\spec D
  \setminus \{0\} > 0$. Then\footnote{The spectral gap condition is
    only need in order to assure that the ranges are closed. If we
    replace $\ran \de$ by $\clo{\ran \de}$ and similarly for $\de^*$,
    we can drop this condition.}
  \begin{gather*}
    \HS   = \ker D \oplus \ran \de^* \oplus \ran \de,\\
    \HS_0 = \ker \de \oplus \ran \de^* \und
    \HS_1 = \ker \de^* \oplus \ran \de.
  \end{gather*}
\end{lemma}
\begin{proof}
  It is a general fact that $\HS_0 = \ker \de \oplus \clo{\ran \de^*}$
  and similarly for $\HS_1$. It remains to show that $\ran \de$ and
  $\ran \de^*$ are closed. Let $\wt D$ be the restriction of $D$ onto
  $(\ker D)^\orth$. By our assumption, $\wt D$ has a bounded inverse,
  namely
  \begin{equation*}
    \wt D^{-1} \cong
    \begin{pmatrix}
      0 & \wt \de^{-1} \\ (\wt \de^*)^{-1} & 0
    \end{pmatrix},
  \end{equation*}
  where $\wt \de$ and $\wt \de^*$ are the restrictions of $\de$ and
  $\de^*$ to $(\ker \de)^\orth$ and $(\ker \de^*)^\orth$,
  respectively. In particular, $\wt \de^{-1}$ and $(\wt \de^*)^{-1}$
  are bounded.

  Let $g \in \clo{\ran \de}$, then there exists a sequence $\{f_n\}_n
  \subset \HS_0$ such that $\de f_n \to g$ in $\HS_1$. Without loss of
  generality, we may assume that $f_n \in (\ker \de)^\orth$.
  Therefore, $\wt \de^{-1} \de f_n = f_n \to \wt \de^{-1} g =: f$.
  Now, $f_n \to f$, $\de f_n \to g$ and $\de$ is closed, so $f \in
  \dom \de$ and in particular, $\de f = g \in \ran \de$.
\end{proof}

\begin{definition}
\label{def:index}    
If~$\ker \de$ and~$\ker \de^*$ are both finite dimensional (i.e., $0
\notin \spec[ess] D$), we define the \emph{index} of~$D$ as
  \begin{equation*}
    \ind D := \dim \ker \de - \dim \ker \de^*.
  \end{equation*}
\end{definition}
Note that~$\ind D$ is the usual Fredholm index of the operator~$\de$.

We need the following fact in order to calculate the index in concrete
examples:
\begin{lemma}
  \label{lem:index}
  Assume that~$\{D_t\}_{t \in \R}$ is a family of bounded Dirac
  operators such that~$t \mapsto D_t$ is norm-continuous. Then~$\ind
  D_t$ is constant.
\end{lemma}
\begin{proof}
  This follows from the fact that the Fredholm index depends
  continuously on the operator and that a continuous function
  into~$\Z$ is locally constant (see
  e.g.~\cite[Lem.~1.4.3]{gilkey:95}).
\end{proof}

We need the notion of a morphism of this structure.
\begin{definition}
  \label{def:susy.mor}
  Suppose that $\HS=\HS_0 \oplus\HS_1$ with operator $\de$ and $\wt \HS
  = \wt \HS_0 \oplus \wt \HS_1$ with operator $\wt \de$ and associated
  Dirac operators $D$ and $\wt D$, respectively, have supersymmetry.
  We say that a linear map $\map \Phi {\dom D}{\dom {\wt D}}$
  \emph{respects supersymmetry} iff $\Phi$ decomposes into
  $\Phi=\Phi_0 \oplus \Phi_1$ where $\Phi_p$ maps $p$-forms onto
  $p$-forms.
\end{definition}

In some cases we need to enlarge the Hilbert space $\HS$ by a space
$\mc N$ on which the exterior derivative acts trivially:
\begin{definition}
  \label{def:dirac.enlarged}
  Let $\mc N$ be a Hilbert space. We set $\HS_{\mc N} := \HS_0 \oplus
  \mc N \oplus \HS_1$. Assume that $\de$ is an exterior derivative on
  $\HS=\HS_0 \oplus \HS_1$. Then we call
  \begin{equation*}
    \map{\de_{\mc N^0}=\de \oplusmerge 0} {\dom \de \oplus \mc N}
    {\HS_1},
    \qquad f \oplus h \mapsto \de f
  \end{equation*}
  the exterior derivative \emph{trivially $0$-enlarged by $\mc
    N$}. The associated Dirac operator will be denoted by $D_{\mc
    N^0}$.

  Similarly, we call
  \begin{equation*}
    \map{\de_{\mc N^1}=\de \oplussplit 0} {\dom \de}
    {\HS_1  \oplus \mc N},
    \qquad f \mapsto \de f \oplus 0
  \end{equation*}
  the exterior derivative \emph{trivially $1$-enlarged by $\mc
    N$}. The associated Dirac operator will be denoted by $D_{\mc
    N^1}$.
\end{definition}
Note that $\de_{\mc N^0}^*=(\de \oplusmerge 0)^* = \de^* \oplussplit
0$ and $\de_{\mc N^1}^*=(\de \oplussplit 0)^*= \de^* \oplusmerge 0$.
Furthermore, $\ker \de_{\mc N^0} = \ker \de \oplus \mc N$, $\ker
\de_{\mc N^0}^* = \ker \de^*$ and $\ker \de_{\mc N^1} = \ker \de$,
$\ker \de_{\mc N^1}^* = \ker \de^* \oplus \mc N$. In particular, we have
\begin{equation}
  \label{eq:dirac.enlarged}
  \ind D_{\mc N^0} = \ind D + \dim \mc N \und
  \ind D_{\mc N^1} = \ind D - \dim \mc N.
\end{equation}
%
\section{Vertex spaces on discrete graphs}
\label{sec:vx.sp}
%

\subsection{Discrete graphs}
\label{sec:dg}
Suppose~$\Graph$ is a discrete weighted graph given
by~$(V,E,\bd,\ell)$ where~$(V,E,\bd)$ is a usual graph, i.e.,~$V$
denotes the set of vertices,~$E$ denotes the set of edges,~$\map \bd E
{V \times V}$ associates to each edge~$e$ the pair~$(\bd_-e,\bd_+e)$
of its initial and terminal point (and therefore an orientation).
That~$\Graph$ is an \emph{(edge-)weighted} graph means that there is a
\emph{length} or \emph{(inverse) edge weight function}~$\map \ell E
{(0,\infty)}$ associating to each edge $e$ a length~$\ell_e$. For
simplicity, we consider \emph{internal} edges only, i.e., edges of
\emph{finite} length~$\ell_e < \infty$.

For each vertex~$v \in V$ we set
\begin{equation*}
  E_v^\pm := \set {e \in E} {\bd_\pm e = v} \qquad \text{and} \qquad
  E_v := E_v^+ \dcup E_v^-,
\end{equation*}
i.e.,~$E_v^\pm$ consists of all edges starting ($-$) resp.\ ending
($+$) at~$v$ and~$E_v$ their \emph{disjoint} union. Note that the
\emph{disjoint} union is necessary in order to allow self-loops, i.e.,
edges having the same initial and terminal point.  The
\emph{(in/out-)}degree of~$v \in V$ is defined as
\begin{equation*}
  \deg^+ v := |E^+_v|, \qquad
  \deg^- v := |E^-_v|, \qquad
  \deg v := |E_v|=\deg^+ v +\deg^- v,
\end{equation*}
respectively. In order to avoid trivial cases, we assume that~$\deg v
\ge 1$, i.e., no vertex is isolated.  On the vertices, we usually
consider the canonical (vertex-)weight $\deg v$ (see e.g.\ the norm
definition of $\lsqr V$ in~\eqref{eq:norm.deg}).

We say that the graph $\Graph$ is \emph{$d$-regular}, iff $\deg v =d$
for all $v \in V$. Furthermore, $\Graph$ is \emph{bipartite}, if there
is a decomposition $V = V_- \dcup V_+$ such that no vertex in $V_-$ is
joined with a vertex in $V_-$ by an edge and similar for $V_+$.

We have the following equalities
\begin{equation}
  \label{eq:vx.ed.bij}
    \bigdcup_{v \in V} E^+_v 
  = \bigdcup_{v \in V} E^-_v 
  = E \und
  \bigdcup_{v \in V} E_v = E \dcup E,
\end{equation}
since each (internal) edge has exactly one terminal vertex and one
initial vertex. In addition, a self-loop edge~$e$ is counted twice in
$E_v$.  In particular,
\begin{equation}
  \label{eq:2nd.graph}
  \sum_{v \in V} \deg v = 2 |E|.
\end{equation}

\subsection{General vertex spaces}
\label{sec:vx.spaces}
We want to introduce a vertex space allowing us to define Laplace-like
operators coming from general vertex boundary conditions for quantum
graphs. The usual discrete Laplacian is defined on~$0$-forms and
$1$-forms, namely, on sections in the trivial bundles
\begin{equation*}
  \Forms^0 \Graph = V \times \C \und
  \Forms^1 \Graph = E \times \C.
\end{equation*}
In order to allow more general vertex boundary conditions in the
quantum graph case later on, we need to enlarge the space at each
vertex~$v$. We denote~$\Gmax_v := \C^{E_v}$ the \emph{maximal vertex
  space} at the vertex~$v \in V$, i.e., a value~$\ul F(v) \in \Gmax_v$
has~$\deg v$ components, one for each adjacent edge.  A (general)
\emph{vertex space} is a family $\{\mc G_v\}_v$ of subspaces~$\mc G_v$
of~$\Gmax_v$ for each vertex $v$. We can consider a vertex space as a
vector bundle
\begin{equation*}
  \Forms^0 \Graph := \bigdcup_{v \in V} \mc G_v
\end{equation*}
over the discrete base space~$V$ with fibres~$\mc G_v$ of mixed rank
generalising the above setting where~$\mc G_v \cong \C$ at each
vertex. An element of~$\Gmax_v$ will generally be denoted by~$\ul F(v)
= \{F_e(v)\}_{e \in E_v}$. Note that
\begin{equation}
  \label{eq:g.max}
   \Gmax := \bigoplus_v \Gmax_v \cong \bigoplus_{e \in E} \C^2
\end{equation}
since each edge occurs twice in the~$E_v$,~$v \in V$ (cf.\
\Eq{vx.ed.bij}).

 We denote by
\begin{equation}
  \label{eq:lsqr.forms}
  \lsqr {\Forms^0 \Graph} = \mc G := \bigoplus_{v \in V} \mc G_v \und
  \lsqr {\Forms^1 \Graph} = \lsqr E 
    = \bigoplus_{e \in E} \frac 1 {\ell_e^{1/2}} \C
\end{equation}
the associated Hilbert spaces of~$0$- and~$1$-forms with norms defined
by
\begin{equation*}
  \normsqr[\mc G] F
  := \sum_{v \in V} |\ul F(v)|^2 
   = \sum_{v \in V} \sum_{e \in E_v}|F_e(v)|^2  \und
  \normsqr[\lsqr E] \eta
  := \sum_{e \in E} |\eta_e|^2\frac 1 {\ell_e}.
\end{equation*}
Abusing the notation, we also call the section space~$\mc G$ a
\emph{vertex space}.

\begin{definition}
  \label{def:local}
  We say that an operator~$A$ on~$\mc G$ is \emph{local} iff~$A$
  decomposes with respect to~$\mc G=\bigoplus_{v \in V} \mc G_v$,
  i.e.,~$A = \bigoplus_{v \in V}A_v$ where~$A_v$ is an operator
  on~$\mc G_v$.
\end{definition}

Associated to a vertex space is an orthogonal projection~$P =
\bigoplus_{v \in V} P_v$ in~$\Gmax$, where~$P_v$ is the orthogonal
projection in~$\Gmax_v$ onto~$\mc G_v$. Alternatively, a vertex space
is characterised by fixing an orthogonal projection~$P$ in~$\mc G$
which is local.
\begin{remark}
  \label{rem:local}
  If $\Graph$ is finite, we can assume without loss of generality that
  $P$ is local. If this is not the case, we can pass to a new graph
  $\wt \Graph$ by identifying vertices $v \in V$ for which $P$ does
  not decompose with respect to $\Gmax_v \oplus \bigoplus_{w \neq v}
  \Gmax_w$. In the worst case, the new graph $\wt \Graph$ is a rose,
  i.e., $\wt \Graph$ consists of only one vertex with $|E|$ self-loops
  attached.
\end{remark}
The following notation will be useful:
\begin{definition}
  \label{def:orient}
  The linear operator~$\map {\tau=\orientmap} \Gmax \Gmax$,~$F
  \mapsto \orient F$, defined by~$\tau := \bigoplus_{v \in V} \tau_v$ and
  \begin{equation*}
    \tau_v (F(v))
    := \orient {\ul F}(v) = \{\orient F_e(v)\}_{e \in E_v}, \qquad
    \orient F_e(v) := \pm F_e(v), \quad \text{if~$v = \bd_\pm e$,}
  \end{equation*}
  is called \emph{orientation map}. We say that~$\tau$ \emph{switches
    from an unoriented evaluation to an oriented evaluation} and vise
  versa.
\end{definition}
Clearly,~$\tau$ is a unitary local involution and given by the
multiplication with~$\orientul \1(v)$ on~$\Gmax_v$ where~$\orient
\1_e(v)=\pm 1$ if~$v=\bd_\pm e$.
\begin{definition}
  \label{def:dual.vx.sp}
  Let~$\mc G= \bigoplus_{v \in V} \mc G_v$ be a vertex space with
  associated projection~$P$. The \emph{dual} vertex space is defined
  by~$\mc G^\orth := \Gmax \ominus \mc G$ with projection~$P^\orth =
  \1 - P$. The \emph{oriented} version of the vertex space~$\mc G$ is
  defined by~$\orient{\mc G} := \tau \mc G$ with projection~$\orient P =
  \tau P \tau$.
\end{definition}
It can easily be seen that $\orient{\mc G}=\mc G$ iff $\orientul
\1(v)=\pm \ul \1(v)$ for all $v \in V$, i.e., iff the graph $\Graph$
is bipartite (with partition $V=V_- \dcup V_+$) and the orientation is
chosen in such a way that $\bd_\pm e \in V_\pm$ for all $e \in E$.

In the following we give several examples of vertex spaces. We will
see later on that these spaces are closely related to quantum graph
Laplacian where the names come from.  We start with two trivial vertex
spaces:
\begin{example}
  \label{ex:vx.sp.triv}

  \indent
  \begin{enumerate}
  \item
    \label{dir}
    We call the trivial subspace~$\mc G_v = \mc G_v^{\min} = 0$ the
    \emph{minimal} or \emph{Dirichlet vertex space}. The corresponding
    projection is~$P_v=0$.
  \item
    \label{neu}
    We call the maximal subspace~$\mc G_v = \Gmax_v$ the
    \emph{maximal} or \emph{Neumann vertex space}. The corresponding
    projection is~$P_v=\1$. Clearly,~$\Gmax$ is dual to~$\mc G^{\min}$.
  \end{enumerate}
\end{example}
These examples are trivial, since every edge decouples from the
others:
\begin{definition}
  \label{def:dec}
  Let~$\mc G_v$ be a vertex space at~$v$ with projection~$P_v$.
  \begin{enumerate}
  \item We say that \emph{$e_1 \in E_v$ interacts with~$e_2 \in E_v$
      in~$\mc G_v$} iff
    \begin{equation*}
      p_{e_1,e_2}(v) :=\iprod {\delta_{e_1}(v)}{P_v
        \delta_{e_2}(v)} \ne 0
    \end{equation*}
    where~$(\delta_{e_1})_e(v) = 1$ if~$e=e_1$ and~$0$ otherwise.
    If~$p_{e_1,e_2}(v) =0$, we say that \emph{$e_1, e_2 \in E_v$
      decouple in~$\mc G_v$.}
  \item We say that~$\mc G_v$ \emph{decouples along~$E_1 \dcup E_2
      \subset E_v$} iff~$e_1$ and~$e_2$ decouple in~$\mc G_v$ for
    all~$e_1 \in E_1$ and~$e_2 \in E_2$.
  \item We say that~$\mc G_v$ is \emph{completely interacting}
    iff~$e_1$ and~$e_2$ are interacting for any~$e_1, e_2 \in
    E_v$,~$e_1 \ne e_2$.
  \end{enumerate}
\end{definition}
\begin{lemma}
  \label{lem:dec}
  The edges~$e_1,e_2 \in E_v$ ($e_1 \ne e_2$) are interacting (resp.\
  decoupling) in~$\mc G_v$ iff they are in~$\mc G_v^\orth$. In
  particular,~$\mc G_v$ is completely interacting iff~$\mc G_v^\orth$
  is.
\end{lemma}
\begin{proof}
  The claim follows immediately from
  \begin{equation*}
        \bigiprod {\delta_{e_1}}{P_v^\orth \delta_{e_2}}
    = - \bigiprod {\delta_{e_1}}{P_v \delta_{e_2}}
   \end{equation*}
   since~$e_1 \ne e_2$.
\end{proof}
\begin{remark}
  \label{rem:dec}
  Let~$\mc G$ be a vertex space associated to the graph~$\Graph$ such
  that~$\mc G_v$ decouples along~$E_1 \dcup E_2 = E_v$, then~$\mc G_v
  = \mc G_{1,v} \oplus \mc G_{2,v}$. Passing to a new graph~$\wt
  \Graph$ with the same edge set~$E(\wt \Graph)=E(\Graph)$ but
  replacing~$v \in V(\Graph)$ by two vertices~$v_1$,~$v_2$
  with~$E_{v_1}=E_1$ and~$E_{v_2}=v_2$, we obtain a new graph with one
  more vertex. Repeating this procedure, we can always assume that no
  vertex space~$\mc G_v$ decouple along a non-trivial
  decomposition~$E_v=E_1 \dcup E_2$. It would be interesting to
  understand the ``irreducible'' building blocks of this decomposition
  procedure.
\end{remark}

We will define now our main example, since it covers many of
classically defined discrete Laplacians on a graph, as we will see
later on:
\begin{definition}
  \label{def:vx.cont}
  We say that a vertex space~$\mc G_v$ is \emph{(weighted) continuous}
  if~$\dim \mc G_v=1$, i.e.,
  \begin{equation*}
    \mc G_v = \C \ul p(v), \qquad |\ul p(v)|^2 = \deg v,
  \end{equation*}
  and~$\mc G_v$ is completely interacting, i.e.,~$p_e(v) \ne 0$ for
  all~$e \in E_v$ where~$\ul p(v)=\{p_e(v)\}_e$.

  A vertex space~$\mc G$ is called \emph{(weighted) continuous} if all
  its components~$\mc G_v$ are (weighted) continuous and if there are
  uniform constants~$p_\pm \in (0,\infty)$ such that
  \begin{equation*}
    p_- \le |p_e(v)| \le p_+, \qquad e \in E_v, \quad V \in V.
  \end{equation*}
  The dual of a continuous vertex space is called an \emph{(unoriented
    weighted) sum vertex space}.
\end{definition}
Applying the procedure of \Rem{dec}, any vertex space~$\mc G_v$ of
dimension~$1$ with generating vector~$\ul p(v)$ has a decomposition of
$\mc G_v$ along~$E_1 := \set{e \in E_v}{p_e(v)\ne 0}$ and~$E_2 := E_v
\setminus E_1$. The corresponding space~$\mc G_{1,v}$ is now a
continuous vertex space.

In all of the following examples, we can choose~$p_\pm=1$ as uniform
bounds.
\begin{example}
  \label{ex:vx.sp.cont}
  \indent
\begin{enumerate}
  \addtocounter{enumi}{2}
  \item
    \label{cont}
    Choosing~$\ul p(v)=\ul \1(v)$, i.e.,~$\mc G_v := \mc G_v^\stand :=
    \C \ul \1(v)= \C(1, \dots, 1)$, we obtain the \emph{(uniform)
      continuous} or \emph{standard} vertex space denoted by~$\mc
    G_v^\stand$ where all coefficients~$p_e(v)=1$.  The associated
    projection is
    \begin{equation*}
      P_v = \frac 1 {\deg v} \mbE
    \end{equation*}
    where~$\mbE$ denotes the square matrix of rank~$\deg v$ where all
    entries equal~$1$. 

  \item
    \label{or.cont}
    \sloppy
    We also have an \emph{oriented} version of the standard vertex
    space, namely~$\mc G^{\orient \stand}=\C \orient \1$
    where~$\orient \1$ is defined in \Def{orient}. In particular,
    \begin{equation*}
      p_e(v) = \pm 1 \quad \text{if} \quad v=\bd_\pm e.
    \end{equation*}

  \item
    \label{sum}
    We call the dual~$\mc G_v^\Sigma := (\mc G_v^\stand)^\orth =
    \Gmax_v \ominus \C(1, \dots, 1)$ of the continuous vertex space
    the \emph{(unoriented uniform) sum} or \emph{$\Sigma$}-vertex
    space.  Its associated projection is
    \begin{equation*}
      P_v = \1 - \frac 1 {\deg v} \mbE.
    \end{equation*}

  \item
    \label{or.sum}
    The \emph{oriented sum vertex space} is the dual of the oriented
    continuous vertex space, i.e.,~$\orient{\mc G}^{\Sigma} :=
    (\orient{\mc G}^{\stand})^\orth$.

  \item 
    \label{magnetic}
    A more general case of continuous vertex spaces is given by
    vectors~$\ul p(v)$ such that~$|p_e(v)|=1$, we call such continuous
    vertex spaces \emph{magnetic}. An example is giving in the
    following way: Let~$\alpha \in \R^E$ be a function associating to
    each edge~$e$ the \emph{magnetic vector potential}~$\alpha_e \in
    \R$ and set
    \begin{equation*}
      p_e(v)= \e^{-\im \orient \alpha_e(v)/2}
    \end{equation*}
    where~$\orient \alpha_e(v):= \pm \alpha_e$ if~$v=\bd_\pm e$ as in
    \Def{orient}. We call the associated vertex space~$\mc
    G_v^{\magn,\alpha}$ \emph{magnetic}.
  \end{enumerate}
\end{example}
\begin{remark}
  \label{rem:vx.sp}
  \noindent
  \begin{enumerate}
  \item Obviously, for the standard vertex space~$\mc G_v^\stand= \mc
    G_v^{\magn,0}$. Furthermore, the oriented standard vertex space
    $\mc G^{\orient \stand}$ of~\eqref{or.cont} is unitary equivalent
    to a special case of magnetic vertex spaces in~\eqref{magnetic}:
    Choose $\alpha_e=\pi$ for all $e \in E$ then $p_e(\bd_\pm e)= \mp
    \im$, i.e., $\ul p(v) = -\im \orientul \1(v)$ and therefore $\mc
    G^{\orient \stand}= \im \mc G^{\magn,\pi}$.

\item
    Note that any magnetic vertex space occurs in the above way:
    Let~$\hat{\mc G}$ be a magnetic vertex space, then~$\hat
    p_e(v)=\e^{-\im \hat A_e(v)}$ for some~$\hat A=\{\ul{\hat A}(v)\}$
    with~$\ul {\hat A}(v)\in \R^{E_v}$. Let
    \begin{equation*}
      \map{\alpha := \dde \hat A} E \R \qquad \text{i.e.,} \qquad
      \alpha_e = \hat A_e(\bd_+ e) - \hat A_e(\bd_- e),
    \end{equation*}
    (we define~$\dde = \dde^{\max}$ in the next section).
    Let~$A_e(v):= \orient \alpha_e(v)/2$, then~$\dde A = \dde \hat A$,
    i.e.,~$A - \hat A \in \ker \dde$. But the kernel of~$\dde$
    consists of the values~$B$ such
    that~$B_e(\bd_+e)=B_e(\bd_-e)=:\beta_e$ for all~$e \in E$
    where~$\beta \in \R^E$, in particular,
    \begin{equation*}
      A_e(v) = \hat A_e(v) + \beta_e.
    \end{equation*}
    Define a unitary map~$F \mapsto \hat F$,~$\hat F_e(v):= \e^{\im
      \beta_e} F_e(v)$ then~$\hat F \in \hat {\mc G}$ iff~$F \in \mc
    G$ where~$\mc G=\mc G^{\magn,\alpha}$ as defined below. In
    particular,~$\hat {\mc G}$ is unitarily equivalent to~$\mc
    G^{\magn, \alpha}$ for some vector potential~$\alpha \in \R^E$.
  \end{enumerate}
\end{remark}
We want to express continuous vertex spaces with respect to the
standard space $\lsqr V$, the ``classical'' space of~$0$-forms~$\map {\wt
  F} V \C$ with norm defined by
\begin{equation}
  \label{eq:norm.deg}
  \normsqr[\lsqr V] {\wt F}
  := \sum_{v \in V} |\wt F(v)|^2 \deg v.
\end{equation}
In particular, the next lemma shows, that the vertex-weight $\deg v$
is canonical in the sense of~\eqref{vx.cont3}:
\begin{lemma}
  \label{lem:vx.cont}
  Let~$\mc G$ be a continuous vertex space with projection~$P$ and
  denote by~$[p^{-1}]$ the operator
  \begin{equation*}
    \map {\bigl[ p^{-1} \bigr]} \Gmax \Gmax, \qquad
    F \mapsto \wt F = \{ \wt {\ul F}(v) \}_v, \qquad
    \wt F_e (v) = \frac {F_e(v)} {p_e(v)}.
  \end{equation*}
  \begin{enumerate}
  \item
    \label{vx.cont1}
    The multiplication operators~$[p^{-1}]$ and~$[p]=[p^{-1}]^{-1}$
    are bounded on~$\Gmax$
  \item
    \label{vx.cont2}
    We have~$[p^{-1}] (\mc G) = \mc G^\stand$ and~$[p^{-1}] (\mc
    G^\orth) = \mc G^{\Sigma |p|^2}$ where
    \begin{equation*}
      \mc G^{\Sigma |p|^2}
      := \Bigset {\wt F \in \Gmax}
                 { \sum_{e \in E_v} |p_e(v)|^2 \wt F_e(v)=0 \quad
                   \forall v \in V}
    \end{equation*}
    for the dual.
  \item 
    \label{vx.cont3}
    Denote~$\map {\wt U} {\mc G^\stand}{\lsqr V}$ the local operator
    mapping~$\wt {\ul F}(v)=\wt F(v)(1,\dots,1)$ onto~$\wt F(v) \in
    \C$, then~$\wt U$ is unitary. Furthermore,
    \begin{equation*}
      \map U {\mc G} {\lsqr V}, \qquad
      U := \wt U \circ \bigl[p^{-1}\bigr]
    \end{equation*}
    is unitary.
    \item
    \label{vx.cont4}
    The transformed projection~$\map{\wt P:= UP} \Gmax {\lsqr V}$ is
    given by
      \begin{equation*}
        (\wt P_v F)(v) 
        = \frac 1{\deg v} \sum_{e \in E_v} \conj {p_e(v)} F_e(v)
        \in \C
      \end{equation*}
      and no coefficient~$p_e(v)$ vanishes.
  \end{enumerate}
\end{lemma}
\begin{proof}
  \eqref{vx.cont1}~The boundedness follows from the global
  bounds~$p_\pm$ on~$|p_e(v)|$ (cf.~\Def{vx.cont}).
  \eqref{vx.cont2}~$[p^{-1}]$ restricted to~$\Gmax_v$ maps the
  vector~$\ul p(v)$ onto~$(1,\dots,1)$, i.e.,~$\mc G_v$ onto~$\mc
  G^\stand_v$; a vector~$\ul F(v) \in \mc G^\orth$ satisfies~$\sum_{e
    \in E_v} \conj{p_e(v)} F_e(v)=0$, and therefore~$\wt{\ul F}(v) \in
  \mc G^{\Sigma|p|^2}_v$. \eqref{vx.cont3}~We have
  \begin{equation*}
    |\wt {\ul F}(v)|_{\C^{E_v}}^2 
    = |\wt F(v)|^2 |(1, \dots, 1)|_{\C^{E_v}}^2
    = |\wt F(v)|^2 \deg v
  \end{equation*}
  and therefore,~$\wt U$ is unitary. Furthermore,
  \begin{equation*}
    \normsqr[\lsqr V]{\wt F}
    = \sum_{v \in V} |\wt F(v)|^2 \deg v
    = \sum_{v \in V} \sum_{e \in E_v} |\wt F(v) p_e(v)|^2
    = \sum_{v \in V} \sum_{e \in E_v} |F_e(v)|^2
    = \normsqr[\mc G] F
  \end{equation*}
  since~$|\ul p(v)|^2 = \deg v$. The last assertion follows by a
  straightforward calculation.
\end{proof}
Note that the decomposition into $\mc G^{\stand}$ and $\mc
G^{\Sigma|p|^2}$ is no longer orthogonal if~$[p^{-1}]$ is not unitary
(i.e.,~$|p_e(v)| \ne 1$ for some~$e \in E_v$).

The trivial, the uniform continuous and the sum vertex spaces are
obviously invariant under permutation of the edges in~$E_v$. Indeed,
these are the only possibilities for such an invariance:
\begin{lemma}
  \label{lem:invariant}
  A vertex space~$\mc G_v$ is invariant under permutation of the
  coordinates~$e \in E_v$ iff~$\mc G_v$ is either maximal
  ($\Gmax_v=\C^{E_v}$), minimal ($\mc G^{\min}_v=0$), uniform continuous
  ($\mc G^\stand_v=\C(1,\dots,1)$) or the sum vertex space ($\mc
  G^\Sigma_v=\C^{E_v} \ominus \C(1,\dots,1)$).
\end{lemma}
\begin{proof}
  It can be shown, that a square matrix~$P$ of dimension~$d=\deg v$ is
  invariant under the symmetric group~$S_d$ of order~$d$ iff~$P$ has
  the form
  \begin{equation*}
    P = a \1 + b \, \mbE,
  \end{equation*}
  since the only subspaces invariant under~$S_d$ are~$\C(1,\dots,1)$
  and its orthogonal complement, and the representation of~$S_d$ on
  the orthogonal complement is irreducible (see e.g.\ the references
  in~\cite{kuchment:04}). Using the relations~$P=P^*$ and~$P^2=P$ for
  an orthogonal projection, we obtain that~$a$ and~$b$ must be real
  and satisfy the relations~$a^2=a$ and~$2ab+(\deg v) b=b$, from which
  the four cases follow.
\end{proof}

\section{Operators on vertex spaces}
\label{sec:op.vx.sp}

In this section, we define a generalised coboundary operator or
exterior derivative associated to a vertex space. We use this exterior
derivative for the definition of an associated Dirac and Laplace
operator in the supersymmetric setting of \Sec{susy}.

\subsection{Discrete exterior derivatives}
\label{sec:dis.ex.der}
On the maximal vertex space~$\Gmax$, we define a general coboundary
operator or \emph{exterior derivative} as
\begin{equation*}
  \map{\dde=\dde^{\max}}{\Gmax}
            {\lsqr E}, \qquad
    (\dde F)_e := F_e(\bd_+ e) - F_e(\bd_- e),
  \end{equation*}

\begin{definition}
  \label{def:discr.ext.der}
  Let~$\mc G$ be a vertex space of the graph~$\Graph$. The exterior
  derivative on~$\mc G$ is defined as
  \begin{equation*}
    \map{\dde_{\mc G}:= \dde^{\max} \restr{\mc G}}{\mc G}
              {\lsqr E}, \qquad
    (\dde F)_e := F_e(\bd_+ e) - F_e(\bd_- e),
  \end{equation*}
  mapping~$0$-forms onto~$1$-forms.
\end{definition}
We often drop the subscript~$\mc G$ for the vertex space, or use other
intuitive notation in order to indicate the vertex space.

We define a multiplication operator~$[\ell^{-1}]$ on~$\Gmax$ and
$\lsqr E$ by
\begin{equation*}
  \bigl(\bigl[\ell^{-1}\bigr] F \bigr)_e(v) = \frac 1 {\ell_e} F_e(v)
     \und
  \bigl(\bigl[\ell^{-1}\bigr] \eta \bigr)_e = \frac 1 {\ell_e} \eta_e,
\end{equation*}
respectively. Clearly,~$[\ell^{-1}]$ is bounded on both spaces iff
there exists~$\ell_0 > 0$ such that
\begin{equation}
\label{eq:len.bd}
  \ell_e \ge \ell_0, \qquad e \in E.
\end{equation}
On a vertex space~$\mc G \le \Gmax$ with associated projection~$P$, we
can relax the condition slightly, namely, we assume that~$P
[\ell^{-1}]$ is bounded, i.e., that
\begin{equation}
  \label{eq:ass.len}
  \kappa:= \sup_{v \in V} \bigl|P_v [\ell^{-1}]_v\bigr|_v < \infty
\end{equation}
where~$|\cdot|_v$ denotes the operator norm for matrices on~$\C^{E_v}$.
\begin{remark} 
  \indent
  \begin{enumerate}
  \item If~\eqref{eq:len.bd} is fulfilled, then~$\kappa \le 1/\ell_0$.
    In particular, if~$\ell_e=\ell_0$ for all~$e \in E$
    then~$\kappa=1/\ell_0$.
    \item For the (uniform) continuous vertex space~$\mc G^\stand$, we
      have
      \begin{equation*}
        |P_v [\ell^{-1}]_v|_v 
        = \frac 1 {\deg v} \sum_{e \in E_v} \frac 1 {\ell_e}.
      \end{equation*}
    \item If we assume that~\eqref{eq:ass.len} holds for $P$ and
      $P^\orth$, then~\eqref{eq:len.bd} is also fulfilled. For
      simplicity, we assume therefore that~\eqref{eq:len.bd} holds (if
      not stated otherwise).
  \end{enumerate}
\end{remark}
\begin{lemma}
  \label{lem:discr.ext.der}
  Assume~\eqref{eq:ass.len}, then~$\dde$ is norm-bounded by~$\sqrt
  {2\kappa}$. The adjoint
  \begin{equation*}
    \map {\dde^*}{\lsqr E}{\mc G}
  \end{equation*}
  fulfills the same norm bound and is given by
  \begin{equation*}
    (\dde^* \eta)(v) 
    = P_v \Bigl( \Bigl\{ \frac 1 \ell_e \orient \eta_e(v) \Bigr\} \Bigr)
    \in \mc G_v,
  \end{equation*}
  where~$\orient \eta_e(v):= \pm \eta_e$ if~$v=\bd_\pm e$ denotes the
  \emph{oriented} evaluation of~$\eta_e$ at the vertex~$v$.
\end{lemma}
\begin{proof}
  We have
  \begin{align*}
    \normsqr[\lsqr E] {\dde F}
    & =   \sum_{e \in E} \frac 1 {\ell_e}
            \bigl| F_e(\bd_+e) - F_e(\bd_-e) \bigr|^2\\
    &\le 2\sum_{v \in V} 
          \Bigl( \sum_{e \in E^+_v} 
                     \frac 1 {\ell_e} \bigl| F_e(v) \bigr|^2 
               + \sum_{e \in E^-_v} 
                     \frac 1 {\ell_e} \bigl| F_e(v) \bigr|^2
          \Bigr) \\
    &\le  2 \sum_{v \in V} 
             \sum_{e \in E_v} \frac 1 {\ell_e} \bigl| F_e(v) \bigr|^2 \\
    &\le  2 \sum_{v \in V} 
             \bigiprod{ [\ell^{-1}]_v \ul F(v)} {\ul F(v)} \\
    &=    2 \sum_{v \in V} 
             \bigiprod{ [\ell^{-1}]_v \ul F(v)} {P_v \ul F(v)} \\
    & \le 2 \kappa \normsqr[\mc G] F
  \end{align*}
  using~\Eq{vx.ed.bij} and the fact that~$\ul F(v) \in \mc G_v$. For
  the second assertion, we calculate
  \begin{align*}
    \iprod {\dde F} \eta 
    & = \sum_{e \in E} \frac 1 {\ell_e}
            \bigl( \conj F_e(\bd_+e) - \conj F_e(\bd_-e) \bigr) \eta_e\\
    & = \sum_{v \in V} 
          \Bigl( \sum_{e \in E^+_v} 
                     \frac 1 {\ell_e}  \conj F_e(v) \, \eta_e
               - \sum_{e \in E^-_v} 
                     \frac 1 {\ell_e}  \conj F_e(v) \, \eta_e
          \Bigr) \\
    & = \sum_{v \in V} 
          \Bigiprod[\Gmax_v] {P_v F} 
            {\Bigl\{ \frac 1 {\ell_e} \orient \eta_e(v) \Bigr\}_{e \in E_v} }
      = \iprod F {\dde^* \eta}
  \end{align*}
  since~$\ul F(v) \in \mc G_v$, i.e.,~$P_v \ul F(v)=\ul F(v)$.
\end{proof}
\begin{example}
  \label{ex:ext.der.triv}
  \indent
  \begin{enumerate}
    \item For the minimal vertex space, we have~$\dde=0$ and~$\dde^*=0$.
    Obviously, these operators are decoupled, i.e., they do not feel
    any connection information of the graph.
  \item For the maximal vertex space, we have (denoting~$\dde = \dde^{\max}$)
    \begin{equation*}
      (\dde^* \eta)_e(v) 
      = \frac 1 \ell_e \orient \eta_e(v).
    \end{equation*}
    The operator~$\dde=\dde^{\max}$ decomposes as~$\bigoplus_e \dde_e$
    with respect to the decomposition of~$\Gmax$ in \Eq{g.max}
    and~$\lsqr E$ in \Eq{lsqr.forms}. Here,
    \begin{equation*}
      \bigl( \map{\dde_e} {\C^2} \C \bigr) \cong
      \begin{pmatrix}
        1 & -1
      \end{pmatrix}
      \und
      \bigl( \map{\dde_e^*} \C {\C^2} \bigr) \cong
      \frac 1 {\ell_e}
      \begin{pmatrix}
        1 \\ -1
      \end{pmatrix}
    \end{equation*}
    where~$F_e=(F_e(\bd_+e), F_e(\bd_-)) \in \C^2$. Again, the
    operators are decoupled, since any connection information of the
    graph is lost.
  \end{enumerate}
\end{example}

\begin{remark}
  \label{rem:embed.edges}
  We can always embed the edge space $\lsqr E$ into $\Gmax$ using the
  operator
  \begin{equation*}
    \map{\iota}{\lsqr E} {\Gmax}, \qquad
    (\iota \eta)_e(v) := \frac 1 {\sqrt {2 \ell_e}} \eta_e.
  \end{equation*}
  Indeed, $\iota$ is an isometry since
  \begin{equation*}
    \normsqr[\Gmax]{\iota \eta} 
     = \frac 12\sum_{v \in V} \sum_{e \in E_v} 
               \frac 1 {\ell_e} |\eta_e|^2
     = \sum_{e \in E}
               \frac 1 {\ell_e} |\eta_e|^2
     = \normsqr[\lsqr E] \eta
  \end{equation*}
  using \Eq{vx.ed.bij}. Furthermore, the range of $\iota$ in $\Gmax$
  is precisely the kernel of $\dde^{\max}$, i.e.,
  \begin{equation*}
    \iota(\lsqr E) = \ker \dde^{\max}
  \end{equation*}
  as it can be checked easily. Moreover, we can write the adjoint of
  the exterior derivative $\dde=\dde_{\mc G}$ on $\mc G$ with
  projection $P$ as
  \begin{equation*}
    \dde^* 
    = P (\dde^{\max})^* 
    =\sqrt 2 P \orient \1 \iota \bigl[\ell^{-1/2}\bigr].
  \end{equation*}
\end{remark}
We can now calculate the exterior derivative and its adjoint in
several general cases. The proofs are straightforward. We start with
the relation to the dual vertex space:
\begin{lemma}
  \label{lem:discr.ext.dual}
  Let~$\mc G$ be a vertex space with exterior
  derivative~$\dde=\dde_{\mc G}$, then
  \begin{align*}
      \map{\dde_{\mc G} \oplusmerge \dde_{\mc G^\orth} 
      &= \dde^{\max}}
      {\Gmax}{\lsqr E}, &
      F \oplus F^\orth 
      &\mapsto \dde_{\mc G}F + \dde_{\mc G^\orth}F^\orth\\
      \map{\dde_{\mc G}^* \oplussplit \dde_{\mc G^\orth}^* 
      &= (\dde^{\max})^*}
      {\lsqr E}{\Gmax}, &
      \eta &\mapsto 
            \dde_{\mc G}^* \eta \oplus \dde_{\mc G^\orth}^*\eta.
  \end{align*}
  In particular,
  \begin{equation*}
      ((\dde_{\mc G^\orth}^* \eta)_e(v) 
      = \frac 1 {\ell_e} \orient \eta_e(v) - (\dde_{\mc G}^* \eta)_e(v).
  \end{equation*}
\end{lemma}

For a continuous vertex space, it is convenient to use the unitary
transformation from $\mc G$ onto $\lsqr V$ (see
\Lemenum{vx.cont}{vx.cont3}):
\begin{lemma}
  \label{lem:discr.ext.cont}
  For a continuous vertex space, the exterior derivative~$\wt \dde :=
  \dde \circ U^{-1}$ transformed back to~$\lsqr V$ is given as
  \begin{equation*}
      (\wt \dde \wt F)_e 
      = p_e(\bd_+e) \wt F(\bd _+ e) 
      - p_e(\bd_-e) \wt F(\bd _- e)
  \end{equation*}
  and its adjoint~$\wt \dde^* = U \circ \dde^*$ by
  \begin{equation*}
      (\wt \dde^*\eta)(v) =
      \frac 1 {\deg v} \sum_{e \in E_v} 
      \frac {\conj{p_e(v)}}{\ell_e} \orient \eta_e(v).
  \end{equation*}
\end{lemma}

Switching the orientation on or off leads to another class of examples:
\begin{lemma}
  \label{lem:discr.ext.or}
  If~$\hat {\mc G}$ is a vertex space with projection~$\hat P$ and if
  we define the ``unoriented'' exterior derivative~$\hat \dde$ via
    \begin{equation*}
      \map {\hat \dde} {\hat {\mc G}} {\lsqr E}, \qquad
      (\hat \dde F)_e := F_e(\bd_+e) + F_e(\bd_-e),
    \end{equation*}
    then its adjoint is given by
    \begin{equation*}
      (\hat \dde ^* \eta) (v)
       = \hat P_v \Bigl( \Bigl\{ \frac 1 \ell_e \eta_e(v) \Bigr\} \Bigr).
    \end{equation*}
    In addition, if~$\mc G = \tau \hat {\mc G}$ is the vertex space
    with switched orientation, then~$\dde = \hat \dde \circ \tau$
    and~$\dde^* = \tau \circ \hat \dde^*$, i.e., the above
    ``unoriented'' exterior derivative~$\hat \dde$ occurs as an
    exterior derivative in the sense of \Def{discr.ext.der} for the
    vertex space~$\tau \mc G$ with switched orientation.
\end{lemma}

We give now some examples of exterior derivatives on continuous vertex
spaces and their duals:
\begin{example}
  \label{ex:discr.ext.cont}
  \indent
  \begin{enumerate}
  \addtocounter{enumi}{2}
    \item For the standard vertex space~$\mc G^\stand$, the exterior
      derivative and its adjoint are unitarily equivalent to
      \begin{equation*}
        \map{\wt \dde}{\lsqr V}{\lsqr E}, \qquad
        (\wt \dde F)_e = F(\bd_+ e) - F(\bd_- e)
      \end{equation*}
      and
      \begin{equation*}
        (\wt \dde^* \eta)(v) 
        = \frac 1 {\deg v} \sum_{e \in E_v} 
                        \frac 1 {\ell_e} \orient \eta_e(v),
      \end{equation*}
      i.e.,~$\wt \dde$ is the classical coboundary operator and~$\wt
      \dde^*$ its adjoint.

    \item If~$\mc G^{\orient \stand} = \tau \mc G^\stand$ is the
      oriented standard vertex space, then the exterior
      derivative~$\dde$ is unitarily equivalent to
      \begin{equation*}
        \map{\wt \dde}{\lsqr V}{\lsqr E}, \qquad
        (\wt \dde F)_e = F(\bd_+ e) + F(\bd_- e)
      \end{equation*}
      and
      \begin{equation*}
        (\wt \dde^* \eta)(v) 
        = \frac 1 {\deg v} \sum_{e \in E_v} 
                        \frac 1 {\ell_e} \eta_e(v).
      \end{equation*}

    \item For the (unoriented) sum vertex space~$\mc G^\Sigma=(\mc
      G^\stand)^\orth$, we have
      \begin{equation*}
        (\dde^* \eta)_e(v) 
        =  \frac 1 \ell_e \orient \eta_e(v) 
        - \frac 1 {\deg v} \sum_{e' \in E_v}
        \frac 1 \ell_{e'} \orient \eta_{e'}(v)
      \end{equation*}

    \item For the (oriented) sum vertex space~$\mc G^{\orient
        \Sigma}=(\mc G^{\orient \stand})^\orth$, we have
      \begin{equation*}
        (\dde^* \eta)_e(v) 
        = \pm \Bigl(\frac 1 \ell_e \eta_e(v) 
           - \frac 1 {\deg v} \sum_{e' \in E_v}
           \frac 1 \ell_{e'} \eta_{e'}\Bigr)
      \end{equation*}
      if~$v = \bd_\pm e$.

    \item For the magnetic vertex space~$\mc G^{\magn, \alpha}$, we
      have
      \begin{equation*}
        \map {\wt \dde}{\lsqr V} {\lsqr E}, \qquad
        (\wt \dde F)_e 
           = \e^{-\im \alpha_e/2} F(\bd_+e)
           - \e^{ \im \alpha_e/2} F(\bd_-e)
      \end{equation*}
      and
      \begin{equation*}
        \map {\wt \dde^*}{\lsqr E}{\lsqr V}, \qquad
        (\wt \dde^* \eta)(v) 
           = \frac 1 {\deg v} \sum_{e \in E_v} 
                \frac 1 {\ell_e} \e^{\im \orient \alpha_e/2}
                        \orient \eta_e(v).
      \end{equation*}
  \end{enumerate}
\end{example}

\subsection{Discrete Dirac operators and Laplacians}
\label{sec:dis.dirac.lap}

Let~$\De=\De_{\mc G}$ be the Dirac operator associated to the exterior
derivative~$\dde=\dde_{\mc G}$ on the vertex space~$\mc G$, i.e.,
\begin{equation*}
  \De =
    \begin{pmatrix}
      0 & \dde^*\\ \dde & 0
    \end{pmatrix}
  \quad\text{with respect to} \quad
  \lsqr {\Forms \Graph} := 
  \lsqr {\Forms^0 \Graph} \oplus \lsqr {\Forms^1 \Graph}
  = \mc G \oplus \lsqr E
\end{equation*}
(cf.~\Def{susy}).

\begin{definition}
  \label{def:discr.laplace}
  We define as in the abstract supersymmetric setting the
  \emph{Laplacians} associated to a vertex space~$\mc G$ as
  \begin{equation*}
    \dlaplacian {\Forms \Graph}   := 
    \dlaplacian {\mc G} := \De_{\mc G} ^2, \qquad
    \dlaplacian {\Forms^0 \Graph} :=
    \dlaplacian [0] {\mc G} := \dde_{\mc G}^* \dde_{\mc G} \und
    \dlaplacian {\Forms^1 \Graph}  :=
     \dlaplacian [1]{\mc G} := \dde_{\mc G} \dde_{\mc G}^*.
  \end{equation*}
\end{definition}
In particular, we have
\begin{subequations}
  \begin{gather}
    \label{eq:lap0}
    (\dlaplacian[0]{\mc G} F)(v)
    = P_v \Bigl( \Bigl\{ \frac 1 \ell_e
         \bigl( F_e(v) - F_e(v_e) \bigr) \Bigr\} \Bigr)\\
    \label{eq:lap1}
    (\dlaplacian[1] {\mc G} \eta)_e
    = \Bigl(P_{\bd_+e} \Bigl( \Bigl\{ \frac 1 \ell_{e'}
         \orient \eta_{e'}(\bd_+e) \Bigr\} \Bigr)
      - P_{\bd_-e} \Bigl( \Bigl\{ \frac 1 \ell_{e'}
         \orient \eta_{e'}(\bd_-e) \Bigr\} \Bigr)\Bigr)_e
  \end{gather}
\end{subequations}
where~$v_e$ denotes the opposite vertex of~$v \in E_v$ on~$e$. Here,
we see that the orientation plays no role for the $0$-form Laplacian.

We have a sort of Hodge decomposition (see \Lem{hodge}):
\begin{lemma}
  \label{lem:hodge.dg}
  Assume that $\De$ has a spectral gap at $0$, i.e., that
  $\dist(0,\spec \De \setminus \{0\} > 0$ (e.g., $\Graph$ finite is
  sufficient). Then
  \begin{gather*}
    \lsqr {\Forms \Graph} 
     = \ker \De \oplus \ran \dde^* \oplus \ran \dde, \qquad \text{i.e.,}\\
    \lsqr {\Forms^0 \Graph} = \mc G 
     = \ker \dde \oplus \ran \dde^* \und
    \lsqr {\Forms^1 \Graph} = \lsqr E
    = \ker \dde^* \oplus \ran \dde.
  \end{gather*}
\end{lemma}

Let us start with the Laplacians acting on the trivial vertex spaces:
\begin{example}
  \indent
  \begin{enumerate}
  \item For the minimal vertex space, we have~$\dlaplacian[p]{\mc
      G^0}=0$ for~$p \in \{0,1\}$.
  \item For the maximal vertex space, we have
    ($\dlaplacian[p]{\max}:=\dlaplacian[p]{\Gmax}$)
      \begin{equation*}
        (\dlaplacian[0]{\max} F)_e(v) = 
          \Bigl\{ 
               \frac 1 {\ell_e} \bigl(F_e(v) - F_e(v_e) \bigr) 
          \Bigr\}_{e \in E_v}.
      \end{equation*}
      The operator~$\dlaplacian[0]{\max}$ decomposes as~$\bigoplus_e
      (\dlaplacian[0]{\max})_e$ with respect to the decomposition
      of~$\Gmax$ in \Eq{g.max}, where
      \begin{equation*}
        \bigl( \map{(\dlaplacian[0]{\max})_e} 
               {\C^2} {\C^2} \bigr) 
        \cong \frac 1 {\ell_e}
        \begin{pmatrix}
          1 & -1 \\ -1 & 1
        \end{pmatrix}.
      \end{equation*}
      Similarly,
      \begin{equation*}
        (\dlaplacian[1]{\max} \eta)_e = \frac 2 {\ell_e} \eta_e,
      \end{equation*}
      i.e.,~$\dlaplacian[1]{\max}=2[\ell^{-1}]$ is a
      multiplication operator on~$\lsqr E$.
    \end{enumerate}
\end{example}

\begin{lemma}
  \label{lem:lap}
  \noindent
  \begin{enumerate}
  \item The Laplacian~$\dlaplacian[p] {\mc G}$ on~$p$-forms associated
    to the vertex space~$\mc G$ is a bounded operator with norm
    bounded by~$2 \kappa$.
  \item On~$1$-forms, we have~$\dlaplacian[1]
    {\max}=\dlaplacian[1]{\mc G} + \dlaplacian[1]{\mc G^\orth}$
    or~$\dlaplacian[1]{\mc G^\orth} = 2 \bigl[\ell^{-1}\bigr] -
    \dlaplacian[1]{\mc G}$ on~$\lsqr E$. In particular, if all
    length~$\ell_e=1$, then
      \begin{equation*}
        \dlaplacian[1]{\mc G^\orth} = 2 - \dlaplacian[1]{\mc G} 
           \und
        \spec{\dlaplacian[1]{\mc G^\orth}} = 2 - \spec{\dlaplacian[1]{\mc G}},
      \end{equation*}
      i.e.,~$\lambda \in \spec {\dlaplacian[1]{\mc G^\orth}}$
      iff~$2-\lambda \in \spec {\dlaplacian[1]{\mc G}}$.
    \item Assume that $\ell_e=1$ then  we have the spectral relation
      \begin{equation*}
              \spec {\dlaplacian[0]{\mc G^\orth}} \setminus \{0,2\}
        = 2 -(\spec {\dlaplacian[0]{\mc G      }} \setminus \{0,2\})
      \end{equation*}
      on~$0$-forms, i.e., if~$\lambda \ne 0,2$, then~$\lambda \in
      \spec {\dlaplacian[0]{\mc G^\orth}}$ iff~$2-\lambda \in \spec
      {\dlaplacian[0]{\mc G}}$.
  \end{enumerate}
\end{lemma}
\begin{proof}
  The first assertion follows immediately from \Lem{discr.ext.der}.
  The second is a consequence of \Lem{discr.ext.dual}. The last
  spectral equality follows from the spectral equality for~$1$-forms
  and supersymmetry to pass from~$1$-forms to~$0$-forms
  (cf.~\Lem{susy}).
\end{proof}
In \Lem{ind.dual} we will prove a relation between the kernels, namely
$\ker \dlaplacian[0]{\mc G^\orth} \cong \ker
\dlaplacian[1]{\orient{\mc G}}$.
\begin{lemma}
  \label{lem:lap.cont}
  Let~$\mc G$ be a continuous vertex space,~$\map{\wt \dde}{\lsqr
    V}{\lsqr E}$ the unitarily equivalent exterior derivative as
  defined in \Lem{discr.ext.cont} and~$\wt \dde^*$ its adjoint,
  then~$\wtdlaplacian[0] {\mc G} := \wt \dde^* \wt \dde$
  and~$\dlaplacian[1]{\mc G}$ are given by
  \begin{gather*}
    (\wtdlaplacian[0]{\mc G} F)(v)
    = \frac 1 {\deg v} \sum_{e \in E_v}
       \frac {\conj {p_e(v)}}{\ell_e} 
        \bigl( p_e(v) F(v) - p_e(v_e) F(v_e)\bigr)\\
    (\dlaplacian[1]{\mc G} \wt \eta)_e
    = - \sum_{e \sim e'} \frac {(\conj p_{e'} p_e)(e \cap e')} 
                               {\ell_{e'}\deg (e \cap e')} 
            \orient \eta_{e'}(e)
      + \Bigl( \frac {|p_e(\bd_+ e)|^2}{\deg \bd_+ e} +
               \frac {|p_e(\bd_- e)|^2}{\deg \bd_- e} 
        \Bigr) \frac 1 {\ell_e} \eta_e,
  \end{gather*}
  where~$e' \sim e$ means that~$e'\ne e$ and~$e'$,~$e$ have the
  vertex~$e' \cap e$ in common. Furthermore,~$\orient
  \eta_{e'}(e)=\eta_{e'}$ if the orientation of~$e$,~$e'$ gives an
  orientation of the path formed by~$e$,~$e'$, and~$\orient
  \eta_{e'}(e)=-\eta_{e'}$ otherwise.
\end{lemma}
We have several important special cases of continuous vertex spaces
and their duals:
\begin{example}
  \label{ex:lap.cont}
  \indent
  \begin{enumerate}
  \addtocounter{enumi}{2}
\item For the standard vertex space~$\mc G^{\stand}$, we have the
  standard (weighted) Laplacian~$\dlaplacian[0]{\stand}$ transformed
  to~$\wtdlaplacian[0]{\stand} = \dlaplacian[0] \Graph =
  \dlaplacian[0] {(\Graph,\ell^{-1})}$ on~$\lsqr V$
  and~$\dlaplacian[1]{\stand}$ on~$\lsqr E$, where
  \begin{gather}
    \label{eq:stand.lap}
        (\wtdlaplacian[0]{\stand} F)(v) =
        \frac 1 {\deg v}
        \sum_{e \in E_v} \frac 1 {\ell_e} \bigl(F(v) - F(v_e) \bigr)\\
        \nonumber
        (\dlaplacian[1]{\stand} \eta)_e 
        = - \sum_{e' \sim e} 
          \frac 1 {\ell_{e'} \deg (e' \cap e)} \orient \eta_{e'}(e)
          + \Bigl( \frac 1 {\deg \bd_+ e} 
                 + \frac 1 {\deg \bd_-e} \Bigr) \frac 1 {\ell_e} \eta_e.
    \end{gather}

  \item For the oriented standard space~$\mc G^{\orient \stand}$, we
    have
    \begin{gather*}
        (\wtdlaplacian[0]{\orient \stand} F)(v) = 
        \sum_{e \in E_v} \frac 1 {\ell_e} \bigl(F(v) + F(v_e) \bigr)\\
        (\dlaplacian[1]{\orient\stand} \eta)_e 
        = \sum_{e' \sim e} 
          \frac 1 {\ell_{e'} \deg (e' \cap e)} \eta_{e'}
          + \Bigl( \frac 1 {\deg \bd_+ e} 
                 + \frac 1 {\deg \bd_-e} \Bigr) \frac 1 {\ell_e}
                 \eta_e.
    \end{gather*}
      Note that 
      \begin{equation}
        \label{eq:lap.std.or}
        \wtdlaplacian[0]{\orient \stand} 
        = 2 [L^\Sigma] - \dlaplacian[0]{\stand},
      \end{equation}
      where~$\dlaplacian[0]{\stand}$ is the standard Laplacian of
      Example~\eqref{cont} and $[L^\Sigma]$ is the multiplication
      operator with
      \begin{equation*}
        L^\Sigma(v) 
        := \frac 1 {\deg v} \sum_{e \in E_v} \frac 1 {\ell_e}.
      \end{equation*}

  \item For the (unoriented) sum vertex space~$\mc G^\Sigma$, the dual
    of~$\mc G^{\stand}$, we have
      \begin{gather*}
        (\dlaplacian[0]{\Sigma} F)_e(v) = 
        \frac 1 {\ell_e} \bigl(F_e(v) - F_e(v_e) \bigr) -
        \frac 1 {\deg v} \sum_{e' \in E_v}
             \frac 1 {\ell_{e'}} \bigl(F_{e'}(v) - F_{e'}(v_{e'}) \bigr)\\
      (\dlaplacian[1]{\Sigma} \eta)_e
        =  \sum_{e' \sim e} 
          \frac 1 {\ell_{e'} \deg (e' \cap e)} \orient \eta_{e'}(e)
          - \Bigl( \frac 1 {\deg \bd_+ e} 
                 + \frac 1 {\deg \bd_-e} - 2 \Bigr) \frac 1 {\ell_e} \eta_e.
    \end{gather*}

    \item For the oriented sum vertex space~$\mc G^{\orient \Sigma}$,
      we have
      \begin{gather*}
        (\dlaplacian[0]{\orient \Sigma} F)_e(v) = 
        \frac 1 {\ell_e} \bigl(F_e(v) - F_e(v_e) \bigr) -
        \frac {\orient \1_e(v)}  
              {\deg v} \sum_{e' \in E_v}
             \frac {\orient \1_{e'}(v)} 
                   {\ell_{e'}} 
               \bigl(F_{e'}(v) - F_{e'}(v_{e'}) \bigr)\\
      (\dlaplacian[1]{\orient \Sigma} \eta)_e
        =  - \sum_{e' \sim e} 
          \frac 1 {\ell_{e'} \deg (e' \cap e)} \orient \eta_{e'}(e)
          - \Bigl( \frac 1 {\deg \bd_+ e} 
                 + \frac 1 {\deg \bd_-e} - 2 \Bigr) \frac 1 {\ell_e} \eta_e.
      \end{gather*}

    \item For the magnetic vertex space~$\mc G^{\magn, \alpha}$, we
      have
      \begin{gather*}
        (\wtdlaplacian[0]{\magn, \alpha} F)(v) = 
        \frac 1 {\deg v}
          \sum_{e \in E_v} \frac 1 {\ell_e} \bigl(F(v) -
                    \e^{-\im \orient \alpha_e(v)} F(v_e) \bigr)\\
        (\dlaplacian[1]{\magn, \alpha} \eta)_e 
        = - \sum_{e' \sim e} 
          \frac {\e^{\im \orient \alpha_{e',e}}} 
                {\ell_{e'} \deg (e' \cap e)} \orient \eta_{e'}(e)
          + \Bigl( \frac 1 {\deg \bd_+ e} 
                 + \frac 1 {\deg \bd_-e} \Bigr) \frac 1 {\ell_e}
                 \eta_e,
      \end{gather*}
      where~$\orient \alpha_{e',e}:= (\orient \alpha_{e'} - \orient
      \alpha_e)(e \cap e')$ denotes the oriented flux along~$e'$
      and~$e$.
  \end{enumerate}
\end{example}

\begin{remark}
  \label{rem:line.graph}
  \sloppy
  The $1$-form Laplacian of \Lem{lap.cont} and especially of
  \Exenum{lap.cont}{or.cont} above can be viewed as an operator on the
  line graph. In order to define the line graph, we assume for
  simplicity, that~$\Graph$ has no self-loops and multiple edges, and
  that no edge is isolated (i.e.,~$\deg \bd_+e$ and~$\deg \bd_-e$ are
  not both equal to~$1$). Let~$L(\Graph)$ be the line graph associated
  to the graph~$\Graph$, i.e,~$V(L(\Graph))=E(\Graph)$ and two
  ``vertices'' in the line graph (i.e., edges in the original
  graph)~$e$,~$e'$ are adjacent iff~$e \ne e'$ and~$e \cap e' \ne
  \emptyset$, i.e., if they meet in a common vertex. We have
  \begin{equation*}
    \deg_{L(\Graph)} e 
    = \deg_{\Graph} \bd_+ e + \deg_{\Graph} \bd_- e -2,
  \end{equation*}
  and in particular, if~$\Graph$ is a~$d$-regular graph,
  then~$L(\Graph)$ is~$(2d-2)$-regular.

  The above example of the~$1$-form Laplacian is a line graph
  Laplacian (up to a multiplication operator with the complex edge
  ``weight''
  \begin{equation*}
    \rho_{e,e'} 
    = \frac {(\conj p_{e'} p_e)(e \cap e') \, \orient \1_{e'}(e) \deg e} 
            {\ell_{e'}\deg (e \cap e')}.
  \end{equation*}
  We will now show how~$\dlaplacian{L(\Graph)}$ becomes a Laplacian
  with positive weights.

  If~$\ell_e=1$ for all edges, then the~$1$-form Laplacian is related
  to the~$0$-form Laplacian~$\dlaplacian[0]{(L(\Graph),\rho)}$ on the
  line graph with edge weights
      \begin{equation*}
        \rho_{e,e'}
        = \frac {\deg_{L(\Graph)} e}{\deg_{\Graph} (e \cap e')}
        = \frac {\deg_{\Graph} \bd_+ e + \deg_{\Graph} \bd_- e - 2}
                {\deg_{\Graph} (e \cap e')}
      \end{equation*}
      via
      \begin{equation*}
          \dlaplacian[1]{\orient \stand}
        = [L] - \dlaplacian[0]{(L(\Graph),\rho)}
      \end{equation*}
      where~$[L]$ is the multiplication operator
      on~$\lsqr{V(L(\Graph))}$ with the function
      \begin{equation*}
        L(e)= \Bigl( \frac 1 {\deg \bd_+ e} 
            + \frac 1 {\deg \bd_-e} \Bigr)
            + \sum_{e' \sim e}  \frac 1 {\deg_{\Graph} (e \cap e')} .
      \end{equation*}
      In particular, if~$\Graph$ is~$d$-regular, then~$L=2$. Moreover,
      $\rho_{e,e'}=(2d-2)/d$ and
      \begin{equation}
        \label{eq:lap.line}
          \dlaplacian[1]{\orient \stand}
        = 2 - \frac {(2d-2)} d \dlaplacian[0]{L(\Graph)}.
      \end{equation}
      where
      now,~$\dlaplacian[0]{L(\Graph)}=\dlaplacian[0]{(L(\Graph),1)}$
      is the line graph Laplacian with edge weights set to~$1$.  In
      addition, we can recover a result
      of~\cite{shirai:00,ogurisu:02}\footnote{Shirai and Ogurisu
        actually showed more: If~$\Graph$ is infinite and~$d \ge 3$,
        then~$d/(d-1)$ is contained in the spectrum of the line graph,
        being an eigenvalue of infinite multiplicity. A corresponding
        eigenfunction lies in~$\ker \dlaplacian[1]{\orient \stand} =
        \ker \dde_{\orient \stand}^*$ by~\eqref{eq:lap.line}. For an
        infinite regular graph, one can see that this space is
        infinite-dimensional (see also \Exenum{ind.discr}{or.cont}).},
      namely a spectral relation for the line graph Laplacian and the
      Laplacian on the graph itself,
      \begin{equation*}
        \spec{\dlaplacian[0]{L(\Graph)}} 
             \setminus \Bigl\{\frac d {d-1} \Bigr\}
        = \frac d {2(d-1)} \bigl( \spec{\dlaplacian[0]{\Graph}} 
                     \setminus \{2\} \bigr),
      \end{equation*}
      using supersymmetry,~\eqref{eq:lap.std.or} and~\eqref{eq:lap.line}. 
      In particular, the spectrum of the line graph is always
      contained in the interval~$[0, d/(d-1)]$ and is therefore not
      bipartite (if~$d \ge 3$).
\end{remark}

\begin{remark}
  \label{rem:subdiv}
  There is another interesting example which relates a Dirac operator
  on $\Graph$ to the (standard) Laplacian on the \emph{subdivision
    graph} $S(\Graph)$ defined as follows
  (cf.~\cite{shirai:00,ogurisu:02}).  Again, we assume for simplicity,
  that $\Graph$ has no self-loops and no double edges and that
  $\ell_e=1$. As vertices we set $V(S(\Graph))=V(\Graph) \dcup
  E(\Graph)$, and the edges are given by $\{v,e\}$ if $v \in \bd e$ in
  the original graph (we do not care about the orientation here).  In
  other words, $S(\Graph)$ is obtained from $\Graph$ by introducing a
  new vertex on each edge. The subdivision graph $S(\Graph)$ is always
  bipartite (choose the above decomposition). If $\Graph$ is
  $d$-regular, then $S(\Graph)$ is $(d,2)$-semiregular, i.e., $\deg
  v=d$ for vertices in $v \in V(\Graph) \subset V(S(\Graph))$ and
  $\deg e=2$ for vertices $e \in E(\Graph) \subset V(S(\Graph))$ with
  respect to the bipartite decomposition.

  The standard Laplacian on $S(\Graph)$ is given as
  \begin{align*}
      (\dlaplacian[0]{S(\Graph)}) H(e)
    &= - \frac 12 \bigl( H(\bd_+e) + H(\bd_-e) \bigr)  + H(e)\\
      (\dlaplacian[0]{S(\Graph)}) H(v) &= - \frac 1d \sum_{e \in E_v}
      H(e) + H(v)
  \end{align*}
  for $H \in \lsqr {V(S(\Graph))}$.  In particular,
  \begin{equation*}
    \dlaplacian[0]{S(\Graph)} \cong
    \begin{pmatrix}
      1 & -\wt \dde^*_{\orient \stand}\\
      -\frac 12 \wt \dde_{\orient \stand} & 1
    \end{pmatrix}
    = \1 -
    \begin{pmatrix}
      1 & 0\\ 0 & \frac 12
    \end{pmatrix} \wt \De_{\orient \stand}
  \end{equation*}
  where $\wt \de_{\orient \stand}$ is the (transformed) exterior
  derivative and $\wt \De_{\orient \stand}$ the Dirac operator
  associated to the oriented standard space
  (cf.~\Exenum{discr.ext.cont}{or.cont}). Furthermore,
  \begin{equation*}
    \bigl( \dlaplacian[0]{S(\Graph)} - 1 \bigr)^2 \cong
    \frac 12
    \begin{pmatrix}
      \wtdlaplacian[0] {\mc G^{\orient \stand}} & 0 \\
      0 & \dlaplacian[1] {\mc G^{\orient \stand}}
    \end{pmatrix} =
    \frac 12
    \begin{pmatrix}
      \wtdlaplacian[0] {\mc G^{\orient \stand}} & 0 \\
      0 & \dlaplacian[1] {\mc G^{\orient \stand}}
    \end{pmatrix}
  \end{equation*}
  and (for $\lambda \ne 1)$ we have $\lambda \in
  \spec{\dlaplacian[0]{S(\Graph)}}$ iff $2(\lambda-1)^2 \in
  \spec{\wtdlaplacian[0] {\mc G^{\orient \stand}}}$ by supersymmetry.
  But the latter operator equals $2 - \laplacian[0]{\Graph}$ by
  \Eq{lap.std.or} since we assumed $\ell_e=1$. Therefore
  \begin{equation*}
    \spec {\dlaplacian[0]{S(\Graph)}} \setminus \{1\}
    = \eta^{-1} \bigr( \spec {\dlaplacian[0] \Graph} 
               \setminus \{2\}\bigl)
  \end{equation*}
  where $\eta(\lambda)=2-2(\lambda-1)^2 =2\lambda( 2 - \lambda)$. One
  can show that $1$ is also an eigenvalue with infinite multiplicity
  of the subdivision graph Laplacian. In particular, we recover again
  a result of~\cite{shirai:00}.
\end{remark}

%
\section{Indices for Dirac operators on discrete graphs}
\label{sec:ind.discr}
%

We start this section with a short excursion into cohomology. Assume
that~$\map{\dde=\dde_{\mc G}}{\mc G}{\lsqr E}$ is an exterior
derivative for the vertex space~$\mc G$.
\begin{definition}
  \label{def:cohom}
  We define the \emph{($\lsqrsymb$-)cohomology} of the graph~$\Graph$
  associated to the vertex space~$\mc G$ as
  \begin{equation*}
    H_{\mc G}^0(\Graph, \C) := \ker \dde \und
    H_{\mc G}^1(\Graph, \C) := \ker \dde^* = \lsqr E \ominus \ran \dde.
  \end{equation*}
  We call
  \begin{equation*}
    b_{\mc G}^p(\Graph) := \dim_\C H_{\mc G}^p(\Graph,\C) \und
    \chi_{\mc G}(\Graph) := b_{\mc G}^0(\Graph) - b_{\mc G}^1(\Graph)
  \end{equation*}
  the \emph{$p$-th Betti-number} and \emph{Euler characteristic}
  associated to the vertex space~$\mc G$, respectively.
\end{definition}

From the definition, it follows that
\begin{equation*}
  \chi_{\mc G}(\Graph) = \ind \De_{\mc G}.
\end{equation*}

In order to derive a sort of ``Gau{\ss}-Bonnet''-theorem, we need the
notion of curvature at a vertex for general vertex spaces:
\begin{definition}
  \label{def:curv.vx}
  We define the \emph{curvature} of the vertex space~$\mc G_v$ at the
  vertex~$v \in V$ as
  \begin{equation*}
    \kappa_{\mc G}(v) := \dim \mc G_v - \frac 12 \deg v.
  \end{equation*}
\end{definition}
The reason for the name will become clear in
\Remenum{ind.discr}{gauss-bonnet}.  Note that there are other notions
of curvature, especially for tessellations (see
e.g.~\cite{baues-peyerimhoff:01}).

In order to calculate the Betti-numbers for a vertex space~$\mc G$, we
need some more notation. For simplicity, we assume that $\Graph$ is
connected.  Let~$\Graph'$ be a \emph{spanning tree} of $\Graph$,
i.e,~$\Graph'$ is simply connected and $V(\Graph')=V(\Graph)$. For
each~$e \in P(\Graph):= E(\Graph) \setminus E(\Graph')$, there exists
a unique cycle~$c_e$ (closed path without repetitions) in~$\Graph$
containing~$e$.
\begin{definition}
  \label{def:prime.cycle}
  A \emph{prime cycle} is a cycle~$c_e$ for some~$e \in P(\Graph)$
  associated to a spanning tree~$\Graph'$ of~$\Graph$ as above. A
  cycle~$c$ is said to be \emph{even/odd} if the number of edges
  in~$c$ is even/odd.
\end{definition}

Before calculating the kernel of the Dirac operator in some examples,
we establish a general result on the dual~$\mc G^\orth$ of a vertex
space~$\mc G$. It shows that actually, $\mc G^\orth$ and the
\emph{oriented} version of $\mc G$ are related:
\begin{lemma}
  \label{lem:ind.dual}
  Assume that the global length bound
  \begin{equation}
  \label{eq:len.2bd}
    \ell_0 \le \ell_e \le \ell_+ \qquad \text{for all $\e \in E$}
  \end{equation}
  holds for some constants $0 < \ell_0 \le \ell_+ < \infty$. Then
  \begin{equation*}
    H^0_{\mc G^\orth}(\Graph,\C)
    = \ker \dde_{\mc G^\orth} 
    \cong \ker \dde_{\orient {\mc G}}^*
    = H^1_{\orient{\mc G}}(\Graph,\C)
  \end{equation*}
  are isomorphic.  In particular, if $\Graph$ is finite, then
  \begin{equation*}
    b^0_{\mc G^\orth}(\Graph) = b^1_{\orient{\mc G}}(\Graph), \qquad 
    b^1_{\mc G^\orth}(\Graph) = b^0_{\orient{\mc G}}(\Graph)  \und
    \chi_{\mc G^\orth} (\Graph) = - \chi_{\orient{\mc G}}(\Graph).
  \end{equation*}
  If in addition, $\De$ has a spectral gap at $0$, then $\ran
  \dde_{\mc G^\orth}^* \cong \ran \dde_{\orient{\mc G^\orth}}$.
\end{lemma}
\begin{proof}
  We define $\map \psi {\lsqr E}{\mc G^\orth}$ via $(\psi \eta)(v):=
  P^\orth_v \{ \frac 1 {\ell_e} \eta \}$, then
  \begin{equation*}
    \normsqr[\Gmax] {\psi \eta}
    = \sum_{v \in V} \sum_{e \in E_v} \frac 1 {\ell_e^2} |\eta_e|^2
    = 2 \sum_{e \in E} \frac 1 {\ell_e^2} |\eta_e|^2
    \le \frac 2 {\ell_0} \normsqr[\lsqr E] \eta,
  \end{equation*}
  and therefore $\psi \eta \in \mc G^\orth$. Furthermore, $\psi(\ker
  \dde_{\orient{\mc G}}^*) \subset \ker \dde_{\mc G^\orth}$ since
  $\dde_{\orient{\mc G}}^* \eta=0$ implies that $P_v\{\frac 1 {\ell_e}
  \eta_e\} = 0$ and therefore, $(\psi \eta)_e(v)=\frac 1 {\ell_e}
  \eta_e$. In particular, the latter expression is independent of
  $\bd_\pm e$, so that $\dde_{\mc G^\orth} (\psi \eta)=0$. The other
  inclusion can be shown similarly: Let $F \in \ker \dde_{\mc
    G^\orth}$ and set $\eta_e := \ell_e F_e(v)$ independent of
  $v=\bd_\pm e$. Then $\eta \in \lsqr E$ using the global \emph{upper}
  bound $\ell_e \le \ell_+$. Furthermore, $(\psi \eta)(v)=P^\orth \ul
  F(v))= \ul F(v)$ and $\dde_{\orient{\mc G}}^* \eta= \orient \1 P
  F=0$ since $F \in \mc G^\orth$. The other assertions follow
  from the definitions and the fact that $\orth$ and
  $\orient \cdot$ are involutions. The isomorphism of the ranges
  follows from \Lem{hodge.dg}.
\end{proof}

When writing the index of~$\De$, we implicitly assume that the graph
is finite, i.e., that~$|E| < \infty$.  We calculate the cohomology for
the list of our examples. For simplicity, we assume that $\Graph$ is
finite and connected. In particular, the global length
bound~\eqref{eq:len.2bd}, i.e.,~$0<\ell_0\le \ell_e \le \ell_+ <
\infty$ is fulfilled.

\begin{example}
  \label{ex:ind.discr}
  \indent
  \begin{enumerate}
  \item For the minimal vertex space, we have~$\ker \dde_{\mc G^0}=0$
    and~$\ker \dde_{\mc G^0}^*=\lsqr E$. In particular,~$\ind \De_{\mc
      G^0} = -|E|$.

  \item \sloppy For the maximal vertex space, we see from
    \Lem{ind.dual} that~$\ker \dde_{\Gmax} \cong \ker \dde_{\mc G^0}^*
    = \lsqr E$ and~$\ker \dde_{\Gmax}^* \cong \ker \dde_{\mc G^0} =
    0$.  In particular,~$\ind \De = |E|$.

  \item For the standard vertex space~$\mc G^\stand$ we obtain the
    classical homology groups~$H^p(\Graph,\C)$. The $0$-th
    Betti-number counts the number of components, i.e.,
    $b_\stand^0(\Graph)=1$, and the $1$-st the number of prime cycles.
    It is a classical fact that~$b_\stand^1(\Graph)= |P(E)|=|E|-|V|+
    1$ and therefore
    \begin{equation*}
      \ind \De_{\mc G^\stand}
      = b_{\stand}^0(\Graph) - b_{\stand}^1(\Graph) 
      = |V| - |E|
      = \chi_{\stand}(\Graph)
    \end{equation*}

  \item For the oriented standard vertex space~$\mc G^{\orient
      \stand}$, the~$0$-th Betti number counts the number of bipartite
    components of~$\Graph$, i.e., if~$\Graph$ is connected,
    then~$b_{\orient \stand}^0=1$ if~$\Graph$ is bipartite and~$0$
    otherwise. This can be seen using the characterisation
    that~$\Graph$ is bipartite iff~$\Graph$ contains no odd cycle.
    Note that~$(\wt \dde \wt F)_e =0$ for each edge $e$ in an odd
    cycle~$c$ implies that~$\wt F$ vanishes on each vertex in~$c$.

    The~$1$-st Betti number counts the number of prime cycles
    $|P(E)|$, where one has to subtract $1$ if there is an odd prime
    cycle. But the existence of an odd (prime) cycle is equivalent to
    the fact that $\Graph$ is not bipartite. In particular,$b_{\orient
      \stand}^1 = |E|-|V|+1$ if $\Graph$ is bipartite and $b_{\orient
      \stand}^1 = |E|-|V|$ otherwise. Again, we have $\ind \De_{\mc
      G^{\orient \stand}}=|V|-|E|$.

  \item For the (unoriented) sum vertex space~$\mc G^\Sigma$ we can
    apply \Lem{ind.dual}. In particular,
    \begin{equation*}
        \ker \dde_{\mc G^\Sigma} 
        = \Bigset {F \in \Gmax} 
               {F_e(\bd_+e) = F_e(\bd_-e) =: F_e, \,\,
               \sum_{e \in E_v} F_e=0}
    \end{equation*}
    is isomorphic to $\ker \dde_{\mc G^{\orient \stand}}$, i.e.,
    $b_\Sigma^0= b_{\orient \stand}^1=|E|-|V|+1$ iff $\Graph$ is
    bipartite and $b_\Sigma^0=|E|-|V|$ otherwise.  Furthermore,
    \begin{equation*}
        \ker \dde_{\mc G^\Sigma}^*
        = \Bigset {\eta \in \lsqr E} 
               { \frac 1 {\ell_e} \orient \eta_e(v) \text{ is independent
                 of~$e \in E_v$ for~$v \in V$}}
    \end{equation*}
    is a sort of ``oriented'' continuity condition. In particular,
    $b_\Sigma^1= b_{\orient \stand}^0=1$ iff $\Graph$ is bipartite and
    $b_\Sigma^1=0$ otherwise.  Finally,
    \begin{equation*}
      \ind \De_{\mc G^\Sigma} = - \chi_{\stand}(\Graph)= |E| - |V|.
    \end{equation*}

  \item The oriented sum vertex space~$\mc G^{\orient \Sigma}$ is dual
    to the standard vertex space, i.e, $b_{\orient
      \Sigma}^0=|E|-|V|+1$ and $b_{\orient \Sigma}^1=1$ by
    \Lem{ind.dual}. In particular,
    \begin{equation*}
      \ind \De_{\mc G^{\orient\Sigma}} = - \chi_{\stand}(\Graph)= |E| - |V|.
    \end{equation*}

  \item Assume that $\Graph$ is finite. For the magnetic vertex
    space~$\mc G^{\magn, \alpha}$, we need to define the flux through
    a circuit $c$: If $c = \sum_{i=1}^n p_i e_i \in H_1(X,\Z)$ is a
    cycle represented in the homology group, then we define the
    \emph{flux} of $\alpha$ through $c$ as
    \begin{equation*}
      c \cdot \alpha := \sum_{i=1}^n p_i \orient \alpha_e(e_i),
    \end{equation*}
    where $\orient \alpha_e(e_i)$ is defined in \Lem{lap.cont}. Now,
    $b^0_{\magn,\alpha}=1$ iff $c \cdot \alpha \in 2\pi \Z$ for all
    cycles $c \in H_1(\Graph,\Z)$ and $b^0_{\magn,\alpha}=0$
    otherwise.  In order to calculate $b^1_{\magn,\alpha}$, we note
    that the linear system~$(\wt \dde_{\magn,\alpha}^* \eta)(v)=0$ for
    all $v \in V$ consists of $|E|$ variables $\eta_e$ and $|V|$
    equations, therefore $b^1_{\magn,\alpha} \ge |E| - |V|$. It
    remains to show that the rank of the coefficient matrix can
    increase by $1$ iff $c \cdot \alpha \in 2\pi \Z$ for all $c \in
    H_1(\Graph,\Z)$, i.e., $b^1_{\magn,\alpha}=|E| - |V| + 1$ in this
    case and $b^1_{\magn,\alpha}=|E| - |V|$ otherwise. We do not give
    a formal proof of this fact here, since the result follows by
    abstract arguments of the next theorem. Indeed, we have
    \begin{equation*}
      \ind \De_{\mc G^{\magn,\alpha}} = \chi_{\stand}(\Graph)
       = |V| - |E|.
    \end{equation*}
  \end{enumerate}
\end{example}

These examples suggest the following theorem:
\begin{theorem}
  \label{thm:index.discr}
  Assume that~$\Graph$ is a finite graph, i.e., that~$|E| < \infty$,
  with vertex space~$\mc G$. Then the index of the Dirac
  operator~$\De$ associated to the exterior derivative~$\dde$ as
  defined in \Def{discr.ext.der} is given by
  \begin{equation*}
    \ind \De = \dim \mc G - |E|.
  \end{equation*}
\end{theorem}
\begin{remark}
  \label{rem:ind.discr}
  \indent
  \begin{enumerate}
  \item
    \label{gauss-bonnet}
    We can interprete the above theorem as a discrete
    ``Gau{\ss}-Bonnet''-theorem for general vertex spaces, namely
    \begin{equation}
      \label{eq:gen.gauss.bonnet}
      \chi_{\mc G} (\Graph)
      = \sum_{v \in V} \kappa_{\mc G}(v)
    \end{equation}
    using \Eq{2nd.graph}, where~$\chi_{\mc G}(\Graph)$ is defined in
    \Def{cohom} and~$\kappa_{\mc G}(v)$ in \Def{curv.vx}.

  \item The index~$\ind \De$ gives at least some simple information on
    the vertex space, namely~$\mc G$ is trivial (i.e.,~$\mc G$ is the
    maximal or minimal vertex space~$\mc G^{\max}$ or~$\mc G^0=0$)
    iff~$\ind \De = \pm |E|$. This follows from \Eq{2nd.graph} and
    \Exenums{ind.discr}{dir}{neu}.

      \item For continuous vertex spaces $\mc G$, we obtain the
        classical case where $\dim \mc G=|V|$, i.e., the classical
        discrete Gau{\ss}-Bonnet formula \Eq{gauss.bonnet}.
  \end{enumerate}
\end{remark}

Before proving our index theorem, we use a deformation argument in
order to calculate the index:
\begin{lemma}
  \label{lem:ind.dg}
  Let $\mc G_0$ and $\mc G_1$ be two vertex spaces with $\dim \mc
  G_0=\dim \mc G_1=n$, then $\ind \De_{\mc G_0}= \ind \De_{\mc G_1}$
  for the Dirac operators associated to the vertex spaces.
\end{lemma}
\begin{proof}
  Denote $P_t$ the associated orthogonal projections, $t \in \{0,1\}$.
  Note that $P_0$ and $P_1$ can be connected by a (norm-)continuous
  path~$P_t$ inside the space of orthogonal projections of rank $n$:
  This can be seen as follows: Let~$\{\phi_{0,k}\}_k$
  and~$\{\phi_{1,k}\}_k$ be two orthonormal bases such that the
  first~$n$ vectors span the range of~$P_0$ and~$P_1$, respectively.
  Let $U_1$ be the unitary operator mapping $\phi_{0,k}$ onto
  $\phi_{1,k}$. Since the space of unitary operators is connected, we
  can find a (norm)-continuous path $t \mapsto U_t$ from the identity
  operator $U_0:=\1$ to $U_1$ such that all operators $U_t$ are
  unitary. Define $P_t:=U_t^* P_0 U_t$, then $t \mapsto P_t$ is a
  continuous path from $P_0$ to $P_1$.

  Let $\phi_{t,k}:= U_t \phi_{0,k}$, then $t \mapsto \{\phi_{t,k}\}_k$
  is a continuous family of orthogonal bases. Let $\De_t$ be the Dirac
  operator defined with respect to $\mc G_t = \ran P_t$. Passing to
  the basis $\{\phi_{t,k}\}_k$ in $\mc G_t$, we may assume that the
  family $t \mapsto \De_t$ is defined on $\C^n \oplus \lsqr E$.
  Moreover, the family $t \mapsto \De_t$ is continuous, since $t
  \mapsto U_t$ is. The index formula follows from \Lem{index}.
\end{proof}

\begin{proof}[Proof of \Thm{index.discr}]

  For each dimension~$n=\dim \mc G$, we use a simple vertex
  space~$\wt{\mc G}$ of dimension~$n$. In particular, we choose the
  space~$\wt {\mc G}_v:= \C^{d_v} \oplus 0 \subset \C^{\deg
    v}=\Gmax_v$ where~$d_v := \dim \mc G_v$ is the dimension of the
  original vertex space at~$v$. This vertex space corresponds to~$d_v$
  ``Neumann'' boundary conditions at the first~$d_v$ edges, and~$\deg
  v-d_v$ ``Dirichlet'' boundary conditions at the remaining edges
  in~$E_v$. Due to the stability of the index shown in \Lem{ind.dg},
  it suffices to calculate the index of the model vertex space~$\wt
  {\mc G}$.

  It is easily seen as in \Exenums{ind.discr}{dir}{neu} that for a
  graph consisting of a single edge with two adjacent vertices, we
  have
  \begin{equation*}
    \ind \De^{\Dir\Dir} = -1, \qquad
    \ind \De^{\Dir\Neu} = 0 \und
    \ind \De^{\Neu\Neu} = 1.
  \end{equation*}
  For example, for the mixed case,~$\dde F= F(1)$ if the Neumann space
  is at the vertex~$1$. In particular,~$\dde$ and~$\dde^*$ are both
  injective, so~$\ker \dde=0$ and~$\ker \dde^*=0$.

  Due to the additivity of the index with respect to orthogonal sums,
  we therefore have
  \begin{equation*}
    \ind \wt \De = \#\{\text{edges with~$\Neu$-$\Neu$}\}
                 - \# \{\text{edges with~$\Dir$-$\Dir$}\}.
  \end{equation*}
  for the Dirac operator associated with the vertex space~$\wt{\mc
    G}$. It remains to show that
  \begin{equation*}
    \label{eq:ind.count}
    \#\{\text{edges with~$\Neu$-$\Neu$}\}
                 - \# \{\text{edges with~$\Dir$-$\Dir$}\}
    = \#\{\text{all~$\Neu$}\} - |E|.
  \end{equation*}
  In order to show this last equality, we argue by induction over the
  dimension~$n$ of~$\wt {\mc G}$, i.e., the total number of Neumann
  conditions in~$\wt {\mc G}$.  For~$n=0$, the vertex space is the
  minimal or Dirichlet vertex space, for which the index formula is
  correct by \Exenum{ind.discr}{dir}. For the induction step~$n \to n+1$
  we have to distinguish two cases:
  \paragraph{Case A} In an existing edge with two Dirichlet or one
  Dirichlet and one Neumann boundary space, we replace one Dirichlet
  space by a Neumann one. This increases the LHS of \Eq{ind.count} as
  well as the RHS by~$1$.
  \paragraph{Case B} We add an edge with Dirichlet and Neumann vertex
  space to the graph. The LHS is unchanged, and in the RHS, we
  increase the number of Neumann conditions by~$1$, but subtract also
  one additional edge.
\end{proof}

%
\section{Exterior derivatives on quantum graphs}
\label{sec:quantum}
%

In this section, we develop the notion of exterior derivatives on
metric graphs. We first start with the definition of a metric graph,
and some general results needed later.

\subsection{Continuous metric graphs}
\label{sec:cont.mg}

A (continuous) metric graph~$\Graph=(V, E, \bd, \ell)$ is formally
given by the same data as a discrete (edge-)weighted graph. The
difference is the interpretation of the space~$X$: We define~$\Graph$
as
\begin{equation*}
  \Graph := \bigdcup_{e \in E} [0,\ell_e] /\sim_\psi
\end{equation*}
where we identify~$x \sim_\psi y$ iff~$\psi(x)=\psi(y)$ with
\begin{equation*}
  \map \psi {\bigdcup_{e \in E} \{0,\ell_e\}} V, 
       \qquad 
  0_e \mapsto \bd_-e, \quad
  \ell_e \mapsto \bd_+e.
\end{equation*}
In the sequel, we often identify the edge~$e$ with the interval
$(0,\ell_e)$ and use~$x=x_e$ as coordinate.  In addition, we denote
$\dd x=\dd x_e$ the Lebesgue measure on~$e$ inducing a natural measure
on~$\Graph$.  The space~$\Graph$ becomes a metric space by defining
the distance of two points to be the length of the shortest path in
$\Graph$ joining these points.

We first define several Hilbert spaces associated with~$\Graph$.
Our basic Hilbert space is
\begin{equation}
  \label{eq:lsqr.qg}
  \Lsqr \Graph := \bigoplussqr_{e \in E} \Lsqr e
\end{equation}
where again~$e$ is identified with~$(0,\ell_e)$.
\begin{remark}
  \label{rem:dg.vs.qg}
  The interpretation of an edge $e$ as a ``continuous''
  \emph{interval} is in contrast with the \emph{discrete} case where
  $e$ is considered as a single point, e.g.~in~\eqref{eq:lsqr.forms}.
  Another point of view is that we use different types of measures; in
  the discrete case a point measure and in the metric case the
  Lebesgue measure. This fact of choosing two different types of
  measures (or even combinations of them) is pointed out in the works
  of Friedman and Tillich
  (cf.~\cite{friedman-tillich:pre04,friedman-tillich:04}) and also
  in~\cite{baker-faber:06,baker-rumely:07}.
\end{remark}

More generally, we define the \emph{decoupled} Sobolev space of order
$k$ by
\begin{equation*}
  \Sobx [k] {\max} \Graph := \bigoplussqr_{e \in E} \Sob[k] e.
\end{equation*}
Obviously, for~$k=0$, there is no difference between~$\Lsqr \Graph$
and the decoupled space. Namely, evaluation of a function at a point
only makes sense if~$k \ge 1$ due to the next lemma. We need the
following notation: For~$f \in \Sobx {\max} \Graph$, we denote
\begin{equation*}
  \ul f = \{ \ul f(v)\}_{v \in V}, \qquad
  \ul f(v) = \{ f_e(v) \}_{e \in E_v}, \qquad
  f_e(v) :=
  \begin{cases}
    f_e(0),      & v = \bd_-e\\
    f_e(\ell_e), & v = \bd_+e
  \end{cases}
\end{equation*}
the \emph{unoriented} evaluation at the vertex~$v$.  Similarly, for~$g
\in \Sobx {\max} \Graph$, we denote
\begin{equation}
  \label{eq:sign}
  \orul g   = \{ \orul g (v)\}_{v \in V}, \qquad
  \orul g (v) = \{ \orient g_e(v) \}_{e \in E_v}, \qquad
  \orient g_e(v) :=
  \begin{cases}
    -g_e(0),      & v = \bd_-e\\
    g_e(\ell_e),  & v = \bd_+e
  \end{cases}
\end{equation}
the \emph{oriented} evaluation at the vertex~$v$.

The following lemma is a simple consequence of a standard estimate for
Sobolev spaces:
\begin{lemma}
  \label{lem:bd.map}
  Assume the condition~\eqref{eq:len.bd} on the edge lengths, i.e.,
  there is $\ell_0 \le 1$ such that $\ell_e \ge \ell_0 >0$ for all $e
  \in E$. Then the evaluation maps
  \begin{equation*}
    \map {(\underline \cdot)} {\Sobx {\max} \Graph} \Gmax, \quad
       f \mapsto \underline f \und
    \map {(\orul \cdot)} {\Sobx {\max} \Graph} \Gmax, \quad
       g \mapsto \orul g,
  \end{equation*}
are bounded by~$2/\sqrt{\ell_0}$.
\end{lemma}
\begin{proof}
  By density, we can assume that $f$ is smooth on each edge. For $e
  \in E_v$, let $\chi_{v,e}$ be the affine linear function with
  value $1$ at $v$ and $0$ at the other vertex $v_e$. Then
  \begin{equation}
    \label{eq:int}
    f_e(v) = 
    \int_e (f_e \chi_{v,e})'(x) \dd x \orient \1_e(v)\\
    = \int_e (f_e' \chi_{v,e})(x) \dd x \orient \1_e(v) 
    + \frac 1 {\ell_e} 
      \int_e f_e(x) \dd x.
  \end{equation}
  In order to avoid an upper bound on $\ell_e$, we replace the edge $e$
  by the shortened edge $\wt e_v$ of length $\wt \ell_e = \max
  \{\ell_e, 1\}$ starting at $v$. Then
  \begin{equation*}
    |f_e(v)|^2 
    \le 2 \wt \ell_e \normsqr[\wt e_v] {f'}  
        + \frac 2 {\wt \ell_e} \normsqr[\wt e_v] f
    \le 2 \max \Bigl\{ 1, \frac 1{\ell_0} \Bigr\} \normsqr[\Sob e] f
  \end{equation*}
  using Cauchy-Schwarz.  Summing the contributions over $e \in E_v$
  and $v \in V$ and using~\eqref{eq:vx.ed.bij} we are done.  The same
  arguments hold for~$\orul g$.
\end{proof}

For a general vertex space~$\mc G$, i.e., a closed subspace of~$\Gmax
:= \bigoplus_{v \in V} \C^{E_v}$, we set
\begin{equation*}
  \Sobx {\mc G} \Graph 
  := \bigset {f \in \Sobx {\max} \Graph} {\ul f \in \mc G} 
  = (\ul \cdot)^- \mc G,
\end{equation*}
i.e., the preimage of~$\mc G$ under the (unoriented) evaluation map,
and similarly,
\begin{equation*}
  \Sobx {\orient {\mc G}} \Graph 
  := \bigset {g \in \Sobx {\max} \Graph} {\orul g \in \mc G} 
  = (\orul \cdot)^- \mc G
\end{equation*}
the preimage of~$\mc G$ under the (oriented) evaluation map. In
particular, both spaces are closed in~$\Sobx {\max} \Graph$.

The reason for two different vertex evaluations becomes clear through
the following lemma:
\begin{lemma}
  \label{lem:part.int}
  For~$f, g \in \Sobx {\max} \Graph$, we have
  \begin{equation*}
      \iprod[\Graph] {f'} g
    = \iprod[\Graph] f {-g'} +
      \iprod[\Gmax] {\ul f} {\orul g}.
  \end{equation*}
\end{lemma}
\begin{proof}
  We have
  \begin{multline*}
      \iprod[\Graph] {f'} g
    + \iprod[\Graph] f {g'}
    = \sum_{e \in E}
        \bigl( \iprod[e]{f'} g + \iprod[e] f {g'}\bigr)\\
    = \sum_{e \in E}
        \bigl[\conj f g\bigr]_{\bd e}
    = \sum_{v \in V} \sum_{e \in E_v}
         \conj f_e(v) \orient g_e(v)
    = \sum_{v \in V} \iprod[\C^{E_v}] {\ul f(v)} {\orul g(v)}
    = \iprod[\Gmax] {\ul f} {\orul g}
  \end{multline*}
  and the latter expression is defined due to \Lem{bd.map}.
\end{proof}

\subsection{Quantum graphs}
\label{sec:qg}
\begin{definition}
  \label{def:laplacian}
  A \emph{Laplacian} on a (continuous) metric graph~$X=(V,E,\bd,\ell)$
  is an operator~$\laplacian X$ acting as~$(\laplacian X f)_e =
  -f_e''$ on each edge~$e \in E$.
\end{definition}

We have the following characterisation
from~\cite[Thm.~17]{kuchment:04}:
\begin{theorem}
  \label{thm:sa}
  Assume the condition~\eqref{eq:len.bd} on the edge lengths, namely
  $\ell_e \ge \ell_0 >0$. Let~$\mc G \le \Gmax$ be a
  (closed) vertex space with orthogonal projection~$P$, and let~$L$ be
  a self-adjoint, bounded operator on~$\mc G$. Then the
  Laplacian~$\laplacian {(\Graph,\mc G,L)}$ with domain
  \begin{equation*}
    \Sob[2]{\Graph,\mc G, L}
    := \bigset{f \in \Sobx[2] {\max}\Graph}
        {\ul f \in \mc G, \quad P \orul f' + L \ul f = 0}
  \end{equation*}
  is self-adjoint with associated quadratic form
  \begin{equation*}
    \qf d_{(\Graph,\mc G,L)}(f)
    := \normsqr[X]{f'} + \iprod[\mc G] f {Lf}
    = \sum_{e \in E} \normsqr[e]{f'_e} 
    + \sum_{v \in V} \iprod[\mc G_v]{ \ul f(v)}{L(v)\ul f(v)}
  \end{equation*}
  and domain~$\dom \qf d_{(\Graph,\mc G,L)}= \Sobx {\mc G} \Graph = \set{f \in
    \Sobx {\max} \Graph}{\ul f \in \mc G}$.
\end{theorem}
\begin{remark}
  \noindent
  \begin{enumerate}
  \item We have a similar assertion for the ``oriented'' version,
    namely, when we replace~$\ul f$ by~$\orul g$ and~$\orul f'$
    by~$\ul g'$.  We will refer to this Laplacian
    as~$\laplacian{(X,\orient {\mc G}, L)}$.
  \item At least for finite graphs, the converse statement is true,
    i.e., if $\Delta$ is a self-adjoint Laplacian in the sense of
    \Def{laplacian} then $\Delta=\laplacian{(\Graph,\mc G, L)}$ for
    some vertex space $\mc G$ and a bounded operator $L$.  For
    infinite graphs, the operator~$L$ may become unbounded but we do
    not consider this case here.
  \item Note that~$\laplacian{(\Graph,\mc G, L)} \ge 0$ iff~$L \ge 0$.
  \end{enumerate}
\end{remark}

We slightly restrict ourselves and consider only those self-adjoint
Laplacians on $\Graph$ that are obtained as in the above theorem:
\begin{definition}
  \label{def:qg}
  A \emph{quantum graph}~$\Graph$ is a metric graph together with a
  self-adjoint Laplacian~$\laplacian {(\Graph,\mc G,L)}$ where~$\mc G$
  (or~$\orient{\mc G}$ for the oriented version) is a vertex space
  and~$L$ a self-adjoint, bounded operator on~$\mc G$. The quantum
  graph is therefore given by~$\Graph=(V,E,\bd,\ell,\mc G, L)$ or by a
  metric graph~$X=(V,E,\bd,\ell)$ and the data~$(X,\mc G, L)$
  (resp.~$(X,\orient{\mc G}, L)$).
\end{definition}

\begin{remark}
  \label{rem:ks}
  In~\cite{kostrykin-schrader:99} (see also~\cite{kps:pre07}) there is
  another way of parametrising all self-adjoint vertex boundary
  conditions, namely for bounded operators $A,B$ on $\Gmax$,
  \begin{equation*}
    \dom \laplacian {(A,B)} =
    \set{ f \in \Sobx[2] \max \Graph} { A \ul f + B \orul f'=0}
  \end{equation*}
  is the domain of a self-adjoint operator $\laplacian {(A,B)}$
  iff
  \begin{enumerate}
  \item $\map{A \oplusmerge B} {\Gmax \oplus \Gmax} {\Gmax}$, $F \oplus
    \orient F \mapsto AF + B \orient F$, is surjective
  \item $AB^*$ is self-adjoint, i.e., $AB^*=BA^*$.
  \end{enumerate}
  Given a vertex space $\mc G \le \Gmax$ and a bounded operator $L$ on
  $\mc G$, we have $\laplacian{(A,B)}=\laplacian{(X,\mc G, L)}$ if we
  choose
  \begin{equation*}
    A \cong
    \begin{pmatrix}
      L & 0\\ 0 & 0
    \end{pmatrix} 
    \und
    B = P \cong
    \begin{pmatrix}
      \1 & 0\\ 0 & 0
    \end{pmatrix}
  \end{equation*}
  with respect to the decomposition $\Gmax = \mc G \oplus \mc
  G^\orth$.  The associated scattering matrix with spectral parameter
  $\mu = \sqrt \lambda$ is
  \begin{equation*}
    S(\mu) 
    = -(A + \im \mu B)^{-1} (A - \im \mu B)
    \cong
    \begin{pmatrix}
      -(L + \im \mu \1)^{-1} (L - \im \mu \1) & 0 \\
      0 & -\1
    \end{pmatrix}.
  \end{equation*}
  In particular, $S(\mu)$ is independent if $\mu$ iff $L=0$, and in
  this case, we have $S(\mu)=\1 \oplus -\1$ for all $\mu$.
\end{remark}
The aim of the subsequent section is to express~$\laplacian
{(\Graph,\mc G,L)}$ as~$\de^* \de$ or~$\de \de^*$. Of course, to do
so, we need~$L \ge 0$ (since operators~$\de^* \de$ and~$\de \de^*$ are
always non-negative). Furthermore, for non-trivial~$L \ne 0$, we need
to enlarge the~$\Lsymb_2$-spaces by the vertex space~$\mc G$.

\subsection{Differential forms, exterior derivatives and Dirac operators}
\label{sec:ext.der}

\begin{definition}
  \label{def:diff.forms}
  For~$p \in \{0,1\}$, let~$\mc G^p$ be a vertex space.
  We call the space
  \begin{equation*}
    \Lsqr {\Forms^p \Graph} := \Lsqr \Graph \oplus \mc G^p
  \end{equation*}
  the \emph{$\Lsymb_2$-space of~$p$-forms}.

  A subspace~$\Sob{\Forms^p \Graph}$ of~$\Lsqr{\Forms^p \Graph}$ is
  called an \emph{$\Sobsymb^1$-space of~$p$-forms}
  iff
  \begin{enumerate}
  \item \label{dense} the space~$\Sob{\Forms^p \Graph}$ is dense
    in~$\Lsqr{\Forms^p \Graph}$,
  \item \label{emb} we have~$\iota^p(\Sobn \Graph) \subset
    \Sob{\Forms^p \Graph}$ and
  \item \label{emb.dual} we have~$(\iota^p)^*(\Sob{\Forms^p \Graph})
    \subset \Sobx {\max} \Graph$,
  \end{enumerate}
  where
  \begin{equation*}
    \embmap{\iota^p}{\Lsqr \Graph} {\Lsqr {\Forms^p \Graph}}, \qquad
    f \mapsto (f,0), \qquad p \in \{0,1\},
  \end{equation*}
  denote the natural embedding operators.  We set
    \begin{equation*}
      \Lsqr {\Forms \Graph} := 
      \Lsqr {\Forms^0 \Graph} \oplus \Lsqr {\Forms^1 \Graph}
          \und
      \Sob{\Forms \Graph} := 
        \Sob{\Forms^0 \Graph} \oplus \Sob{\Forms^1 \Graph}.
    \end{equation*}
\end{definition}
Note that~$(\iota^p)^*$ is the projection onto the first factor.

For the definition of an exterior derivative, we need the following
decoupled operator~$\map{\de_0}{\Sobn \Graph}{\Lsqr{\Graph}}$,
$\de_0f=f'$. Note that~$\de_0$ is a closed operator with adjoint
$\de_1g = -g'$ and~$\dom \de^* = \Sobx {\max}\Graph$. We now define an
exterior derivative associated to spaces of~$p$-forms. Examples are
given below.
\begin{definition}
  \label{def:ext.der}
  Let~$\Lsqr{\Forms^p \Graph}$ be an~$\Lsymb_2$-space of~$p$-forms
  and~$\Sob{\Forms^0 \Graph}$ be an~$\Sobsymb^1$-space of~$0$-forms.
  We call an operator
  \begin{equation*}
    \map \de {\Sob {\Forms^0 \Graph}}{\Lsqr {\Forms^1 \Graph}}
  \end{equation*}
  an \emph{exterior derivative} on the metric graph~$\Graph$ iff the
  following conditions are fulfilled:
  \begin{enumerate}
  \item \label{closed}
    \sloppy The operator~$\de$ is closed as unbounded
    operator from the~$0$-form space~$\Lsqr {\Forms^0 \Graph}$ into
    the~$1$-form space~$\Lsqr {\Forms^0 \Graph}$.
  \item \label{emb.ext} We have~$\de \iota^0 = \iota^1 \de_0$,
    i.e.,~$\de (\iota^0 f)=\iota^1 f'$ for all~$f \in \Sobn \Graph$.
  \item \label{emb.dual.ext} We have~$(\iota^1)^* \de = -\de_0^*
    (\iota^0)^*$, i.e.,~$(\iota^1)^* \de \wt f = ( (\iota^0)^* \wt
    f)'$.
  \end{enumerate}
\end{definition}
Note that the closeness of~$\de$ ensures that we choose the ``right''
norm on~$\Sob{\Forms^p \Graph}$ (and not an artificially smaller
space).
\begin{lemma}
  \label{lem:ext.der}
  Given the~$p$-form spaces~$\Lsqr {\Forms^p \Graph}$,~$p \in
  \{0,1\}$, the~$0$-form space $\Sob{\Forms^0 \Graph}$ and an exterior
  derivative~$\map \de {\Sob{\Forms^0 \Graph}}{\Lsqr{\Forms^0
      \Graph}}$, then the adjoint~$\de^*$ is uniquely defined and
  closed as operator from~$\Lsqr{\Forms^1 \Graph}$
  into~$\Lsqr{\Forms^0 \Graph}$. Its domain
  \begin{equation*}
    \Sob{\Forms^1 \Graph} := \dom \de^*
  \end{equation*}
  is an~$\Sobsymb^1$-space of~$1$-forms (cf.\
  \Defenums{diff.forms}{dense}{emb.dual}).
\end{lemma}
\begin{proof}
  Since~$\de$ is densely defined by \Defenum{ext.der}{dense}, it
  follows, that~$\de^*$ is uniquely determined, closed and densely
  defined, i.e., \Defenum{diff.forms}{closed} is fulfilled. In order
  to verify \Defenum{diff.forms}{emb}, we have to show that for~$f \in
  \Sobn \Graph$, the~$1$-form~$\iota^1 f$ is in~$\dom \de^*$: Set~$h
  := - \iota^0 \de_0 f$. Then~$h \in \Lsqr{\Forms^0 \Graph}$ and we
  have
  \begin{equation*}
    \iprod h g
    = \iprod {-\de_0^* (\iota^0)^* \wt f} f 
    = \iprod {(\iota^1)^* \de \wt f} f
    = \iprod { d \wt f} {\iota^1 f}
  \end{equation*}
  for all~$\wt f \in \dom d$, i.e.,~$\iota^1 f \in \dom \de^*$, where
  we used \Defenum{ext.der}{emb.dual.ext}. Condition~\eqref{emb.dual}
  follows similarly from \Defenum{ext.der}{emb.ext}.
\end{proof}
\begin{definition}
  \label{def:dirac}
  We call the operator
  \begin{equation*}
    \map D {\Sob{\Forms \Graph}}{\Lsqr {\Forms \Graph}}, \qquad
    D(f,g):=(\de^* g, \de f)
  \end{equation*}
  the \emph{Dirac-operator} associated to the~$p$-form
  spaces~$\Sob{\Forms \Graph} \subset \Lsqr {\Forms \Graph}$ and the
  exterior derivative~$\de$.
\end{definition}

\begin{remark}
  \label{rem:data.given}
  \indent
  \begin{enumerate}
    \item Obviously,~$D$ is a closed and self-adjoint operator with the
      matrix representation
      \begin{equation*}
        D \cong
        \begin{pmatrix}
          0 & \de^* \\ \de & 0
        \end{pmatrix}
      \end{equation*}
      where we split the space in its~$0$- and~$1$-form component.
    \item In order to define a Dirac operator~$D$, it suffices --- due
      to \Lem{ext.der} --- to determine
    \begin{equation*}
      \mc G^0, \quad \mc G^1, \quad
      \Sob{\Forms^0 \Graph}, \quad
      \map{\de}{\Sob{\Forms^0 \Graph}} {\Lsqr {\Forms^1 \Graph}}
    \end{equation*}
    and to ensure that the conditions of
    \Defenums{diff.forms}{dense}{emb.dual} for~$p=0$ and
    \Defenums{ext.der}{closed}{emb.dual.ext} are fulfilled. In this
    case, the Dirac operator is uniquely determined.
  \end{enumerate}
\end{remark}

Again, we have a ``baby'' version of the Hodge decomposition theorem
(see \Lem{hodge}):
\begin{lemma}
  \label{lem:hodge.qg}
  Assume that $D$ has a spectral gap at $0$, i.e., that $\dist(0,\spec
  D \setminus \{0\} > 0$ (e.g., $\Graph$ compact is sufficient). Then
  \begin{gather*}
    \Lsqr {\Forms \Graph} 
     = \ker D \oplus \ran \de^* \oplus \ran \de, \qquad \text{i.e.,}\\
    \Lsqr {\Forms^0 \Graph}
     = \ker \de \oplus \ran \de^* \und
    \Lsqr {\Forms^1 \Graph}
    = \ker \de^* \oplus \ran \de.
  \end{gather*}
\end{lemma}

We will now give concrete examples of Dirac operators. Since at this
stage it is not clear what definition is ``natural'' we list some
reasonable possibilities how to define the~$\Sobsymb^1$-spaces:

\begin{lemma}
  \label{lem:dirac}
  \indent
  \begin{enumerate}
  \item
    \label{simple}
    \textbf{The simple case:} We set
    \begin{equation*}
      \mc G^0 := 0, \quad \mc G^1 := 0, \quad
      \Sob{\Forms^0 \Graph} := \Sobx {\mc G} \Graph, \quad
      \de f := f'.
    \end{equation*}
    Then
    \begin{equation*}
      \Sob{\Forms^1 \Graph} = \Sobx {\orient{\mc G}^\orth} \Graph \und
      \de^*g=-g'.
    \end{equation*}
  \item 
    \label{0.enlarged}
    \textbf{The~$0$-enlarged space:} Let~$\mc G$ be a vertex space
    with bounded operator~$L \ge 0$ on~$\mc G$. We set~$\mc G^0 := \mc
    G$,~$\mc G^1 := 0$,
    \begin{equation*}
      \Sob{\Forms^0 \Graph}
      := \bigset{(f,F) \in \Sobx {\mc G} \Graph \oplus \mc G}
         {\ul f = L^{1/2} F}
    \end{equation*}
    and~$\de (f,F) = f'$.  Then we have
    \begin{gather*}
      \Sob {\Forms^1 \Graph} := \Sobx {\max} \Graph 
             \und 
      \de^* g = (-g',L^{1/2}P \orul g).
    \end{gather*}
    
  \item
    \label{0.enlarged.proj}
    \textbf{The~$0$-enlarged space with projection:} Let~$\mc G$
    be a vertex space with associated projection~$P$ and bounded
    operator~$L \ge 0$ on~$\mc G$. We set~$\mc G^0 := \mc G$,~$\mc G^1
    := 0$,
    \begin{equation*}
      \Sob{\Forms^0 \Graph}
      := \bigset {(f,F) \in \Sobx {\max} \Graph \oplus \mc G} 
                 {P \ul f = L^{1/2} F}
    \end{equation*}
    and~$\de (f,F) = f'$.  Then
    \begin{equation*}
      \Sob{\Forms^1 \Graph} = \Sobx {\orient {\mc G}} \Graph
         \und
      \de^* g = (-g', L^{1/2} \orul g).
    \end{equation*}

  \item
    \label{1.enlarged}
    \textbf{The~$1$-enlarged space:} Let~$\mc G$ be a vertex space
    with bounded operator~$L \ge 0$. We set
    \begin{equation*}
       \mc G^0 = 0, \quad  \mc G^1:= \mc G,  \quad
      \Sob {\Forms^0 \Graph} := \Sobx {\max} \Graph,
             \quad
      \de f = (f, L^{1/2} P \ul f).
    \end{equation*}
    Then
    \begin{equation*}
      \Sob{\Forms^1 \Graph} = 
      \bigset{(g,G) \in \Sobx {\max} \Graph \oplus \mc G}
             { \orul g \in \mc G, \quad \orul g + L^{1/2} G =0}
    \end{equation*}
    and~$\de^* (g,G)= -g'$.

  \item
    \label{1.enlarged.proj}
    \textbf{The~$1$-enlarged space with projection:} Let~$\mc G$
    be a vertex space with associated projection~$P$ and bounded
    operator~$L \ge 0$ on~$\mc G$. We set~$\mc G^0 := \mc G$,~$\mc G^1
    := 0$,
    \begin{equation*}
      \Sob {\Forms^0_{\mc G} \Graph} := \Sobx {\mc G} \Graph, \quad 
      \de f = (f',L^{1/2} \ul f).
    \end{equation*}
    Then
    \begin{equation*}
      \Sob{\Forms^1 \Graph} = 
      \bigset{(g,G) \in \Sobx {\max} \Graph \oplus \mc G}
             { P \orul g + L^{1/2} G =0}
    \end{equation*}
    and~$\de^* (g,G)= -g'$.
  \end{enumerate}
\end{lemma}
\begin{proof}
  We only check the conditions for~\eqref{0.enlarged} since the other
  cases are similar. We apply \Lem{ext.der} and have to show first
  that $\Sob{\Forms^0\Graph}$ is an $\Sobsymb^1$-space of $0$-forms
  and second, that $\de$ is an exterior derivative. In order to show
  the first, note that $\Sob{\Forms^0\Graph}$ is dense in
  $\Lsqr{\Forms^0\Graph}$: Let $(f,F) \in \Lsqr \Graph \oplus \mc G$
  and $\eps>0$. By density of $\Sobx {\max} \Graph$ we can find a function
  $f_1 \in \Sobx {\max} \Graph$ such that $\norm[\Lsqr \Graph] {f - f_1}
  \le \eps/2$. Furthermore, we can change $f_1$ to $f_2$ near a vertex
  $v$ in such a way that $\ul f_2(v)=F(v)$ and that their norm
  difference does not exceed $\eps/2$. Then $(f_2,F)$ has distance at
  most $\eps$ from $(f,F)$ in $\Lsqr {\Forms^0 \Graph}$.

  The second and third condition of \Def{diff.forms} are obviously
  fulfilled. In order to show that $\de$ is an exterior derivative we
  have to check the conditions of \Def{ext.der}. For the closeness of
  $\de$ note that the graph norm of $\de$ defined by
  \begin{equation*}
    \normsqr[\de]{(f,F)}
    := \normsqr[\Lsqr{\Forms^1}] {\de f}
     + \normsqr[\Lsqr{\Forms^0 \Graph}] {(f,F)}
    = \normsqr[\Lsqr \Graph] {f'}
     + \normsqr[\Lsqr \Graph] f + \normsqr[\mc G] F
  \end{equation*}
  is the Sobolev norm. It remains to show that $\Sob{\Forms^0\Graph}$
  is closed in $\Sobx {\max} \Graph \oplus \mc G$: Note that $(f,F)
  \mapsto P\ul f - L^{1/2} F$ is continuous by \Lem{bd.map} and since
  $L$ is bounded. Furthermore, $\Sob{\Forms^0 \Graph}$ is the kernel
  of this map and therefore closed.  The second and third condition of
  \Def{ext.der} follow by an immediate calculation.
\end{proof}

\begin{remark}
  \label{rem.l.inv}
  By splitting the vertex space $\mc G$ we may assume that $L$ is
  invertible on a smaller vertex space: For example, in
  the~$0$-enlarged case~\eqref{0.enlarged}, the condition $\ul f =
  L^{1/2} F$ implies that $\ul f \in (\ker L^{1/2})^\orth=: \mc G_1$.
  Then we have
  \begin{equation*}
      \Sob{\Forms^0 \Graph}
      = \bigset {(f,F) \in \Sobx {\mc G_1} \Graph \oplus \mc G_1} 
                 {\ul f = L_1^{1/2} F_1} \oplus \mc G_0
  \end{equation*}
  where $\mc G_0 := \ker L$ and $L_1 := L \restr {\mc G_1}$ is
  invertible.
\end{remark}
\begin{remark}
  \label{rem:l.infty}
  We assume here that $L$ is invertible on $\mc G$. Then we can pass
  to the limit $L \to \infty$ in the equation $G=L F$ in the following
  sense: We consider the limit $L^{-1} \to 0$ in $L^{-1}G = F$, i.e.,
  $F=0$ and no restriction on $G$. We use this interpretation in order
  to show how the above different cases are related in the limit case:
  \begin{enumerate}
    \addtocounter{enumi}{1}
  \item \textbf{The~$0$-enlarged space:} Here, the condition in
    $\Sob{\Forms^0 \Graph}$ is $\ul f = L^{1/2}F$. The limit $L \to
    \infty$ leads to $F=0$, i.e., $\mc G^0=0$. Moreover, the second
    component in $\de^* g$ has to vanish, i.e., $P \orul g=0$.  In
    particular, the added space $\mc G^0= 0$ vanishes and we
    arrive at the simple case~\eqref{simple}.

  \item \textbf{The~$0$-enlarged space with projection:} The condition
    in $\Sob{\Forms^0 \Graph}$ is $P \ul f = L^{1/2}F$. The limit $L
    \to \infty$ leads again to $F=0$, i.e., $\mc G=0$. Furthermore,
    the second component in $\de^* g$ has to vanish, i.e., $\orul
    g=0$.  In particular, we arrive at the simple case~\eqref{simple}
    with Neumann boundary space $\mc G=\Gmax$.

  \item
    \textbf{The~$1$-enlarged space:} The condition in $\Sob{\Forms^1
      \Graph}$ is $\orul g + L^{1/2}G=0$. The limit $L \to \infty$
    here leads to $G=0$, i.e., $\mc G^1=0$ and therefore $P \ul f=0$.
    In particular, we arrive at the simple case~\eqref{simple} with
    the roles of $\mc G$ and $\mc G^\orth$ interchanged.

  \item \textbf{The~$1$-enlarged space with projection:} Finally, in
    this case, we arrive at the simple case~\eqref{simple} with
    Dirichlet vertex space $\mc G=\mc G^{\min}=0$.
  \end{enumerate}
\end{remark}

As in \Sec{susy} we can associate a \emph{Laplacian}
$\laplacian{\Forms \Graph} := D^2$ to the Dirac operator~$D$ on the
metric graph~$X$ with differential form space~$\Sob{\Forms \Graph}
\subset \Lsqr {\Forms \Graph}$ with natural domain
\begin{equation*}
  \Sob[2]{\Forms \Graph} 
   = \bigset {\omega \in \Sob{\Forms \Graph}} 
             {D\omega \in \Sob{\Forms \Graph}}.
\end{equation*}
Furthermore, the operator decomposes into the components
\begin{equation*}
    \laplacian {\Forms \Graph} 
  = \laplacian {\Forms^0 \Graph} \oplus
    \laplacian {\Forms^1 \Graph} 
  = \de^* \de \oplus \de \de^*
\end{equation*}
with natural domains
\begin{gather*}
  \Sob[2]{\Forms^0 \Graph} := \bigset{ f \in \Sob {\Forms^0 \Graph}}
  { \de f \in \Sob{\Forms^1 \Graph}}\\
  \Sob[2]{\Forms^1 \Graph} := \bigset{ g \in \Sob {\Forms^1 \Graph}} {
    \de^* g \in \Sob{\Forms^0 \Graph}}.
\end{gather*}

\begin{lemma}
  \label{lem:laplace}
  Assume that~$D$ is a Dirac operator on the metric graph~$\Graph$ and
  that~$\laplacian{\Forms \Graph} := D^2$ is its associated Laplacian.
  Then the components~$\laplacian {\Forms^p \Graph}$ act as Laplacians
  on the metric graph part, i.e.,
  \begin{gather*}
    \map{(\iota^0)^* \laplacian {\Forms^0 \Graph} \iota^0
      = \de_0^* \de_0} {\Sobx[2] \Dir \Graph}{\Lsqr \Graph}\\
    \map{(\iota^1)^* \laplacian {\Forms^1 \Graph} \iota^1 =
      \de_0^* \de_0} {\Sobx[2] \Dir \Graph}{\Lsqr \Graph},
  \end{gather*}
  where~$\de_0^* \de_0 = \bigoplus_{e \in E} \laplacianD e$ is the sum
  of the Dirichlet Laplacian on each edge.
\end{lemma}
\begin{proof}
  The assertion follows directly from
  \Defenums{ext.der}{emb.ext}{emb.dual.ext}.
\end{proof}

We give the concrete domains of the Laplacians in each of the above cases:
\begin{lemma}
  \label{lem:laplace.ex}
  \indent
  \begin{enumerate}
  \item \textbf{The simple case:} Here we have
  \begin{gather*}
    \Sob[2]{\Forms^0 \Graph} := 
    \bigset{ f \in \Sobx[2] {\max} \Graph} 
           {\ul f \in \mc G, \quad \orul f' \in \mc G^\orth}, \qquad
        \laplacian{\Forms^0 \Graph} f = -f''\\
    \Sob[2]{\Forms^1 \Graph} := 
    \bigset{ g \in \Sobx[2] {\max} \Graph} 
           {\orul g \in \mc G^\orth, \quad \ul g' \in \mc G}, \qquad
        \laplacian{\Forms^1 \Graph} g = -g''.
  \end{gather*}
  In particular, on~$0$-forms, we have the quantum graph~$(X,\mc G,
  0)$, and on~$1$-forms the quantum graph~$(X, \orient{\mc G}^\orth,0)$.
  \item \textbf{The~$0$-enlarged space:} Here we have
    \begin{gather*}
      \Sob[2]{\Forms^0 \Graph} := 
      \bigset{ (f,F) \in \Sobx[2] {\max} \Graph \oplus \mc G} 
             {\ul f \in \mc G, \ul f = L^{1/2} F},\\
      \Sob[2]{\Forms^1 \Graph} := 
         \bigset{g \in \Sobx[2] {\max} \Graph} 
             {\ul g' \in \mc G, \quad \ul g' + LP\orul g = 0},\\
       \laplacian{\Forms^0 \Graph} (f,F) = (-f'',L^{1/2} P \orul f'),
               \qquad
         \laplacian{\Forms^1 \Graph} g = -g''.
     \end{gather*}
     The~$1$-form space~$\Sob[2]{\Forms^1 \Graph}$
     equals~$\Sob[2]{\Graph, \orient{\mc G}^{\max}, \wt L}$ where~$\wt
     L (F,F^\orth):= (LF,0)$ with respect to the decomposition~$\Gmax
     = \mc G \oplus \mc G^\orth$. In particular, the~$1$-form space
     represents the quantum graph~$(\Graph, \orient{\mc G}^{\max}, \wt
     L)$.
  \item \textbf{The~$0$-enlarged space with projection:} We have
    \begin{gather*}
      \Sob[2]{\Forms^0 \Graph}
       := \bigset{(f,F) \in \Sobx[2] {\max} \Graph \oplus \mc G}
                 { P \ul f = L^{1/2} F, \quad \orul f' \in \mc G}\\
      \Sob[2]{\Forms^1 \Graph}
       := \bigset{g \in \Sobx[2] {\max} \Graph}
                 { \orul g \in \mc G, \quad P \ul g' + L \orul g =0},\\
      \laplacian {\Forms^0 \Graph} (f,F) = (-f'',L^{1/2} \orul f')
         \und
      \laplacian {\Forms^1 \Graph} g = -g''.
    \end{gather*}
    The~$1$-form Laplacian defines an (oriented) quantum
    graph~$(\Graph,\orient {\mc G}, L)$.

  \item \textbf{The~$1$-enlarged space:} We have
    \begin{gather*}
      \Sob[2]{\Forms^0 \Graph}
       := \bigset{f \in \Sobx[2] {\max} \Graph}
                 { \orul f' \in \mc G, \quad \orul f' + L P\ul f =0}\\
      \Sob[2]{\Forms^1 \Graph}
       := \bigset{(g,G) \in \Sobx[2] {\max} \Graph \oplus \mc G}
                 { \orul g \in \mc G,  \quad \orul g + L^{1/2} G=0},\\
      \laplacian {\Forms^0 \Graph} f = -f''
         \und
      \laplacian {\Forms^1 \Graph} (g,G) = (-g'',-L^{1/2} P \ul g').
    \end{gather*}
    We have~$\Sob[2]{\Forms^0 \Graph}=\Sob[2]{\Graph, \Gmax, \wt L}$
    where~$\wt L$ is defined as in Case~\eqref{0.enlarged}. In
    particular, the~$0$-form Laplacian defines the quantum
    graph~$(\Graph, \Gmax, \wt L)$.

  \item \textbf{The~$1$-enlarged space with projection:} We have
    \begin{gather*}
      \Sob[2]{\Forms^0 \Graph}
       := \bigset{f \in \Sobx[2] {\max} \Graph}
                 { \ul f \in \mc G, \quad P \orul f' + L \ul f =0}\\
      \Sob[2]{\Forms^1 \Graph}
       := \bigset{(g,G) \in \Sobx[2] {\max} \Graph \oplus \mc G}
                 { \ul g' \in \mc G,  \quad P \orul g + L^{1/2} G=0},\\
      \laplacian {\Forms^0 \Graph} f = -f''
         \und
      \laplacian {\Forms^1 \Graph} (g,G) = (-g'',-L^{1/2} \ul g').
    \end{gather*}
    The~$0$-form Laplacian defines the ``classical'' quantum
    graph~$(\Graph, \mc G, L)$.
  \end{enumerate}
\end{lemma}

%
\section{Index formulas for metric graphs}
\label{sec:kernel}
%
\subsection{Isomorphism between kernels of discrete and quantum graph
  Dirac operators}
\label{sec:iso.ker}
We will now present one of the main results of this article, namely we
establish an isomorphism between~$\ker D$ and~$\ker \De$ in the five
cases of differential forms mentioned above respecting the
supersymmetric space decomposition.

We need some notation: For a bounded operator~$L$ in $\mc G$, we
set~$\mc G_0 := \ker L$ and~$\mc G_1 := \mc G \ominus \ker L$.
Furthermore,~$\mc G_i^\orth := \Gmax \ominus \mc G_i$.  In addition,
we denote the projections corresponding to $\mc G_i$ and $\mc
G_i^\orth$ by $P_i$ and $P_i^\orth$, respectively. Similarly, we
denote the corresponding exterior derivatives by $\map{\dde_i}{\mc G_i}
{\lsqr E}$ and $\map{\dde_i^\orth}{\mc G_i^\orth}{\lsqr E}$.

The discrete Dirac operator needs to be trivially enlarged by the
space $\mc N=\mc G_0$ in the cases
\eqref{0.enlarged}--\eqref{1.enlarged.proj} (cf.\
\Def{dirac.enlarged}).  Finally,~$\eta := \int g$ is defined by
$\eta_e := \int_0^{\ell_e} g_e(x) \dd x$.

\begin{theorem}
  \label{thm:index}
  Assume~\eqref{eq:len.bd}, then
  \begin{equation*}
    \map \Phi {\Sob {\Forms \Graph}}{\lsqr {\Forms \Graph}}
  \end{equation*}is a bounded operator with
  norm bounded by~$2/\sqrt {\ell_0}$, and~$\Phi(\ker D)=\ker \De$ is
  an isomorphism respecting the supersymmetry (i.e., $\Phi=\Phi_0
  \oplus \Phi_1$, cf.~\Def{susy.mor}). In particular, $\Phi_0(\ker
  \de)=\ker \dde$ and $\Phi_1(\ker \de^*) = \ker \dde^*$ are isomorphisms
  and
  \begin{equation*}
    \ind D = \ind \De,
  \end{equation*}
  where~$D$ and $\De$ are the Dirac operators associated to the
  exterior derivatives $\de$ and $\dde$, respectively. Furthermore,
  $\dde$ and $\Phi$ are given in the
  following cases:
  \begin{enumerate}
  \item \textbf{The simple case:} Here, $\map \dde {\mc G} {\lsqr E}$,
  \begin{gather*}
    \map \Phi {\Sobx {\mc G} \Graph \oplus \Sobx {\orient {\mc G}} \Graph}
              {\mc G \oplus \lsqr E}, \qquad
                  (f,g) \mapsto (\ul f, \textstyle \int g),\\
    \ind D = \ind \De = \dim \mc G - |E|.
  \end{gather*}
\item \textbf{The~$0$-enlarged space:} Here, 
    \begin{gather*}
    \map \dde {\mc G_1 \oplus \mc G_0} {\lsqr E}, \qquad
        F_1 \oplus F_0 \mapsto \dde_1 F_1, \\
    \map \Phi {\Sob {\Forms \Graph}}
        {(\mc G_1 \oplus \mc G_0) \oplus \lsqr E}, \qquad
                  (f,F,g) \mapsto (P_1 \ul f, P_0 F, \textstyle \int g),\\
    \ind D = \ind \De = \dim \mc G - |E|.
    \end{gather*}

  \item \textbf{The~$0$-enlarged space with projection:} We have
    \begin{gather*}
      \map \dde {\mc G_0^\orth \oplus \mc G_0} {\lsqr E}, \qquad
         F_0^\orth \oplus F_0 \mapsto \dde_0^\orth F_0^\orth,\\
    \map \Phi {\Sob {\Forms \Graph}}
        {(\mc G_0^\orth \oplus \mc G_0) 
                 \oplus \lsqr E}, \qquad
                  (f,F,g) \mapsto (P_0^\orth \ul f, P_0 F, 
                   \textstyle \int g),\\
    \ind D = \ind \De = |E|.
    \end{gather*}

  \item \textbf{The~$1$-enlarged space:} We have 
    \begin{gather*}
    \map \dde {\mc G_1^\orth} {\mc G_0 \oplus \lsqr E}, \qquad
        F_1^\orth \mapsto 0 \oplus \dde_1^\orth F_1^\orth,\\
    \map \Phi {\Sob {\Forms \Graph}}
        {\mc G_1^\orth \oplus (\mc G_0 \oplus \lsqr E)}, \qquad 
                  (f,g,G) \mapsto (P_1^\orth \ul f, P_0 G, 
                   \textstyle \int g),\\
    \ind D = \ind \De = |E| - \dim \mc G.
    \end{gather*}

  \item \textbf{The~$1$-enlarged space with projection:} We have
    \begin{gather*}
    \map \dde {\mc G_0} {\mc G_0 \oplus \lsqr E}, \qquad
        F_0 \mapsto 0 \oplus \dde_0 F_0,\\
    \map \Phi {\Sob {\Forms \Graph}}
        {\mc G_0 \oplus (\mc G_0 \oplus \lsqr E)}, \qquad
                  (f,g,G) \mapsto (P_0 \ul f, P_0 G, 
                   \textstyle \int g),\\
    \ind D = \ind \De = - |E|.
    \end{gather*}
  \end{enumerate}
\end{theorem}
\begin{remark}
  \indent
  \begin{enumerate}
  \item Note that in all cases, the index is independent of~$L$, i.e.,
    of the decomposition of~$\mc G$ into~$\mc G_0=\ker L$ and~$\mc
    G_1=\mc G \ominus \ker L$ as one expects since the index should be
    constant passing to the limit $L \to 0$.
  \item In the first two cases, we obtain the Euler characteristic as
    index (if~$\mc G= \mc G^{\stand}$). These two cases are the ones
    we obtain by a limit argument where the metric graph is approached
    by a manifold (cf.~\cite{exner-post:05} and a forthcoming paper)
    provided the transversal manifold~$F$ is simply connected (see
    also \ExS{std}{borderline}).
  \item If we assume that~$L$ is invertible, then the index in each
    case remains the same when passing to the limit $L \to \infty$
    (cf.\ \Rem{l.infty}).
 \item We can interprete $\Phi$ in the above theorem as a sort of
   Hilbert chain morphism (cf.~\cite[Ch.~1]{lueck:02}). For example,
   in the $0$-enlarged case~\eqref{0.enlarged}, we have
   \begin{equation*}
     \begin{diagram}
     0 & \rTo   & \Sob {\Forms^0 \Graph}      & \rTo^\de     &
                  \Lsqr{\Forms^1 \Graph} & \rTo & 0\\
           &   & \dTo{\Phi_0}               &              &
                  \dTo{\Phi_1}               &      &     \\
     0 & \rTo   & \mc G_1 \oplus \mc G_0      & \rTo^{\dde_1 \oplusmerge 0} &
                  \lsqr E                & \rTo   & 0
    \end{diagram}
   \end{equation*}
   where the rows are obviously chain complexes (with bounded maps)
   and the diagram is commutative. Note that indeed,
   $\map{\Phi_1}{\Lsqr {\Forms^1\Graph}}{\lsqr E}$, $g \mapsto \int
   g$, is a bounded map also on the $\Lsqrsymb$-space.  The
   commutativity of the diagram follows from the fact that
   \begin{align*}
     \bigr((\dde \Phi_0 - \Phi_1 \de)(f,F)\bigl)_e
     &= (P_1 f)_e(\bd_+ e) - (P_1 f)_e(\bd_- e) - \int_e f_e' \dd x\\ 
     &= - (P_0 f)_e(\bd_+ e) + (P_0 f)_e(\bd_- e)
   \end{align*}
   for $f \in \Sob {\Forms^0 \Graph}$. But note that $\ul f =
   L^{1/2}F$ implies $\ul f \in (\ker L)^\orth=\mc G_1$, so that $\dde
   \Phi_0 = \Phi_1 \de$, i.e., $\Phi$ is a chain morphism. The
   corresponding homology induces the above isomorphism of the Dirac
   operator kernels. The other cases can be treated similarly. We will
   stress this abstract point of view (and also an interpretation of
   the ``enlarged'' spaces as twisted chain complexes) in a
   forthcoming publication.
 \end{enumerate}
\end{remark}

\begin{proof}
  The boundedness of~$\Phi$ in the particular cases follows
  immediately from \Lem{bd.map} and Cauchy-Schwarz. We only prove the
  second case, since the other ones are similar. First, we note that
  in Case~\eqref{0.enlarged}, we have~$(f,F,g) \in \ker D$
  iff~$f_e$,~$g_e$ are constant,~$\ul f \in \mc G_1$,~$\orul g \in \mc
  G_1^\orth$ and~$\ul f = L^{1/2} F$. In particular,~$\dde \ul
  f=0$,~$\dde^* \int g = P_1 \orul g = 0$ and therefore~$\Phi(\ker D)
  \subset \ker \De$.

  In order to show that~$\Phi$ is injective on~$\ker D$, we note
  that~$P_1 \ul f=0$ implies~$\ul f=0$ (since already~$\ul f \in \mc
  G_1$). Furthermore,~$P_0 F=0$ and from~$0=\ul f=L^{1/2} F$ we
  conclude~$P_1 F=0$, i.e.,~$F=0$. Finally,~$\int g=0$
  and~$g_e=\const$ implies that~$g=0$.

  It remains to show that~$\ker \De \subset \Phi(\ker D)$,
  let~$(F_1,F_0,\eta) \in \ker \De$. Set~$f_e(x):= F_{1,e}(\bd_\pm e)$
  (both values are the same since~$\de F_1=0$),~$F:= L^{-1/2} F_1 +
  F_0$ and~$g_e(x):=\eta_e$. Then~$(f,F,g) \in \Sob{\Forms \Graph}$
  since~$\ul f \in \mc G_1 \le \mc G$ and~$L^{1/2} F = F_1 = \ul f$.
  Next, $D(f,F,g)=(-g',L^{1/2} P \orient g,f')=0$ since~$L^{1/2} P_0
  \orient g=0$ and~$L^{1/2} P_1 \orient g=L^{1/2} \de^* \eta=0$,
  and~$f_e$,~$g_e$ are constant.  Finally,~$\Phi(f,F,g)=(P_e \ul f,
  P_0 F, \int g)=(F_1,F_0,\eta)$ and the assertion is proven.

  The index formulas follow from \Thm{index.discr} and
  \Eqs{dirac.enlarged}{2nd.graph}.
\end{proof}

\begin{remark}
  \label{rem:curv.qg}
  We would like to interprete the above index formula together with
  the discrete index formula of \Thm{index.discr} as a
  ``Gau{\ss}-Bonnet theorem'' on quantum graphs. To do so, we define
  the \emph{curvature} $\kappa:=\kappa_{(\Graph,\mc G,L)}$ of the
  quantum graph $(\Graph, \mc G, L)$ as
  \begin{equation}
    \label{eq:curv.qg}
    \kappa_e(x) 
    = 2   \Bigl( 
        \frac {\kappa_{\wt{\mc G}} (\bd_-e)}   
             {\sum_{e' \in E_{\bd_-e}} \ell_{e'}}  (\ell_e - x) +
        \frac {\kappa_{\wt{\mc G}} (\bd_+e)} 
             {\sum_{e' \in E_{\bd_+e}} \ell_{e'}}  x \Bigr)
  \end{equation}
  where $\kappa_{\wt {\mc G}}(v)$ is the \emph{discrete} curvature
  defined in \Def{curv.vx}, but for the vertex space $\wt {\mc G}$
  given by $\wt {\mc G}_v = \mc G_v$ in the
  cases~\eqref{simple}--\eqref{0.enlarged}, by $\wt {\mc G}_v =
\Gmax_v$ in the case~\eqref{0.enlarged.proj}, by $\wt {\mc G}_v = \mc
G_v^\orth$ in the case~\eqref{1.enlarged}, and by $\wt {\mc G}_v = 0$
in the last case~\eqref{1.enlarged.proj} (cf.\ the index formulas in
\Thm{index}).  In particular, we can interprete the index formula
$\ind D = \ind \De$ as
  \begin{equation}
    \label{eq:gauss.bonnet.qg}
    \ind D = \int_{\Graph} \kappa \dd x
  \end{equation}
  since 
  \begin{equation}
    \label{eq:curv2}
    \int_{\Graph} \kappa \dd x =
    \sum_{v \in V} \kappa_{\wt{\mc G}}(v) = \ind \De
  \end{equation}
  by an obvious calculation and \Thm{index.discr}. Note that the
  choice of $\kappa_e$ is somehow arbitrary, but it is the unique way
  to define it if we require that~\eqref{eq:curv2} holds, that
  $\kappa_e(v)=c(v) \kappa_{\wt G}(v)$ for a sequence $c(v)>0$ and
  that $\kappa_e''=0$.

  In particular, we have
  \begin{equation*}
    \kappa(v) = \kappa_{(\Graph,\mc G,L)}(v)
       = \frac {2 \kappa_{\wt {\mc G}} (v)} {\sum_{e' \in E_v} \ell_{e'}}.
  \end{equation*}
  If we have a continuous vertex space, i.e., $\dim \mc G_v=1$ like
  the standard vertex space, then $\kappa(v)=0$ iff $\deg v=2$. This
  reflects the fact that a vertex of degree $2$ is invisible.
  Furthermore if $\deg v=1$ (i.e., a ``dead end'' with Neumann
  boundary space), then $\kappa(v)>0$. Furthermore, if $\deg v \ge 3$,
  then $\kappa(v) < 0$. Moreover, shorter lengths $\ell_e$ at a vertex
  $v$ mean a higher absolute value of the curvature.

  For example, a dead end $e$ with Dirichlet boundary space at $v \in
  \bd e$ has negative curvature. In some sense, one could say that
  high negative curvature forces the function to vanish: If the dead
  end has length $\ell_e \to 0$ with standard vertex conditions on the
  other vertex $w \in \bd e$ (of degree $\ge 3$), then $e$ has
  curvature $\kappa_e \to -\infty$ as $\ell \to 0$, and finally forces
  the function to vanish also on $w$.

  On the other hand a dead end $e$ of length $\ell_e \to 0$ with
  Neumann boundary space at the endpoint $v$ has curvature tending to
  $\infty$, but the curvature at the other point $w$ is negative and
  remains finite. Therefore $\kappa_e(x)=0$ for a point $x \to w$ as
  $\ell \to \infty$, and here, the dead end just ``disappears'' in the
  limit.
\end{remark}

%
%

\subsection{Metric graphs as limits of smooth spaces}
\label{sec:mg.smooth}
We will give several examples of quantum graph operators on
$X_0=\Graph$ which occur as limits of an appropriate smooth
approximation $X_\eps$. A simple example is given if $X_0$ is embedded
in $\R^2$ and if we choose some open neighbourhood $X_\eps$ of $X_0$.
Note that
\begin{equation*}
  \chi(X_\eps) = \chi(X_0),
\end{equation*}
since $X_0$ and $X_\eps$ are homotopy-equivalent.

In~\cite{rubinstein-schatzman:01,kuchment-zeng:01,kuchment-zeng:03%
  ,exner-post:05,post:05,post:06}, the convergence of the $0$-form
operators has been established in various situations. We will show in
a forthcoming article, that the result extends also to differential
forms on $X_\eps$ under suitable conditions. Note that in the three
first examples below, the ``approximating'' Laplacian $\de_\eps^*
\de_\eps$ on $X_\eps$ (with Neumann boundary conditions on $\bd
X_\eps$) and its dual $\de_\eps \de_\eps^*$ on $1$-forms have index
equal to $\chi(X_0)$ (more precisely, the Dirac operator associated to
the exterior derivative $\map{\de_\eps}{\Sob{X_\eps}}{\Lsqr{\Forms^1
    X_\eps}}$ has index equal to $\chi(X_0)$). 

We indicate the limit operators acting on a metric graph in several
situations:
\begin{example}[Standard boundary conditions]
  \label{ex:std}
  If the vertex neighbourhoods do not shrink too slow (e.g., the
  $\eps$-neighbourhood of the embedded metric graph $X_0 \subset \R^2$
  is good enough), then the Neumann Laplacian on functions converges
  to the standard metric graph Laplacian (see the references above),
  and we will also show that the $1$-form Laplacian on $X_\eps$
  converges to the $1$-form metric graph Laplacian. In particular, the
  vertex space of the limit operator is $\mc G=\mc G^{\stand}$ and the
  domains are given by
  \begin{align*}
    \Sobx [2] \stand \Graph 
    &:= \bigset{ f \in \Sobx[2] {\max} \Graph}
    {\text{$\ul f(v)$ independent of $e \in E_v$, 
                 $\sum_{e \in E_v} \orul f'(v)=0$}},\\
    \Sobx [2] {\orient \Sigma} \Graph 
    &:= \bigset{ g \in \Sobx[2] {\max} \Graph}
    {\text{$\ul g'(v)$ independent of $e \in E_v$,
                 $\sum_{e \in E_v} \orul g(v)=0$}}
  \end{align*}
  as domains for the Laplacian on $0$- and $1$-forms, respectively.
  The associated Dirac operator $D$ has index equal to the Euler
  characteristic~(see~\eqref{eq:kernel.std}). The same is true for
  more general Schr\"odinger operators on $X_\eps$ like magnetic
  Laplacians (for the convergence, see
  e.g.~\cite{kuchment-zeng:01,exner-post:07}).  Magnetic Laplacians
  have been studies throughoutly in~\cite{kostrykin-schrader:03}.
\end{example}

\begin{example}[The decoupling case]
  \label{ex:dec.dir}
  In~\cite{kuchment-zeng:03} and~\cite[Sec.~6]{exner-post:05} there is
  a class of approximations $X_\eps \subset \R^2$ (roughly with slowly
  decaying vertex neighbourhood volumes of order $\eps^{2\alpha}$ with
  $0<\alpha < 1/2$). In this case, the limit operator on $0$-forms is
  \begin{equation*}
    \bigoplus_{e \in E} \laplacianD e \oplus
    \bigoplus_{v \in V} 0,
  \end{equation*}
  i.e., the $0$-enlarged case~\eqref{0.enlarged} with $\mc G=\mc
  G^{\stand}$ and operator $L=0$. Again, the index of the associated
  Dirac operator is $\chi(\Graph)$
  (cf.~\Thmenum{index}{0.enlarged}). The dual operator is the
  decoupled Neumann operator.
\end{example}

\begin{example}[The borderline case]
  \label{ex:borderline}
  In~\cite{kuchment-zeng:03} and~\cite[Sec.~7]{exner-post:05} there is
  a special class of approximations $X_\eps \subset \R^2$ where the
  volume of a vertex neighbourhood $U_\vxeps$ is $\vol U_\vxeps= \eps
  \vol U_v$ (i.e., $\alpha=1/2$). In this case, the limit operator on
  $0$-forms is of the form \Lemenum{laplace.ex}{0.enlarged} with $\mc
  G=\mc G^\stand$ and $L(v)=(\vol U_v)^{-1}$ (multiplication
  operator). In particular, the ``bizzar'' boundary conditions in this
  case with the enlarged graph space are ``natural'' in this setting.
  Again, the index of the associated Dirac operator is $\chi(\Graph)$.
  Note that the dual operator is a ``real'' quantum graph Laplacian,
  namely the domain consists of functions $g \in \Sobx [2]{\max}
  \Graph$ such that
  \begin{equation*}
    g' \quad \text{is continuous}, \qquad
    (\deg v)(\vol U_v) \ul g'(v) = \sum_{e \in E_v} \orul g(v),
  \end{equation*}
  i.e., a type of $\delta'$-condition with strength given by the local
  volume (and with oriented evaluation, since we are on $1$-forms).
\end{example}

\begin{example}[The Dirichlet decoupling case]
\label{ex:dir}
In~\cite{post:05} we proved an approximation result for the Laplacian
with Dirichlet boundary conditions on a certain set $X_\eps \subset
\R^2$ which is ``small'' around the vertex neighbourhood. The limit
operator on functions in this case is the simple decoupled operator
$\bigoplus_e \laplacianD e$, i.e., the simple case~\eqref{simple} with
$\mc G=0$. The index formula in this case leads to
  \begin{equation*}
    \ind D = - |E|.
  \end{equation*}
  Note that the index of the Dirichlet Laplacian on $X_\eps$ is the
  \emph{relative Euler characteristic}
  \begin{equation*}
    \chi(X_\eps, \bd X_\eps)= \chi(X_\eps)-\chi(\bd X_\eps)
    = \chi(X_\eps) = \chi(X_0)=|V|-|E|
  \end{equation*}
  in this case, which indicates that the $1$-form Laplacian on
  $X_\eps$ in this case does not converge to the $1$-form Laplacian
  $\laplacian[1]{\mc G^{\min}}=\bigoplus_e \laplacianN e$.  We will
  treat this question also in a forthcoming publication.
\end{example}


\begin{thebibliography}{FKW07}

\bibitem[BF06]{baker-faber:06}
M.~Baker and X.~Faber, \emph{Metrized graphs, {L}aplacian operators, and
  electrical networks}, Quantum graphs and their applications, Contemp. Math.,
  vol. 415, Amer. Math. Soc., Providence, RI, 2006, pp.~15--33. 

\bibitem[BP01]{baues-peyerimhoff:01}
O.~Baues and N.~Peyerimhoff, \emph{Curvature and geometry of tessellating plane
  graphs}, Discrete Comput. Geom. \textbf{25} (2001), no.~1, 141--159.
  

\bibitem[BR07]{baker-rumely:07}
M.~Baker and R.~Rumely, \emph{Harmonic analysis on metrized graphs}, Canad. J.
  Math. \textbf{59} (2007), no.~2, 225--275. 

\bibitem[CdV98]{colin:98} Y. Colin~de Verdi{\`e}re, \emph{Spectres de
    graphes}, Cours Sp\'ecialis\'es, vol.~4, Soci\'et\'e
  Math\'ematique de France, Paris, 1998. \MR{99k:05108}

\bibitem[Chu97]{chung:97}
Fan R.~K. Chung, \emph{Spectral graph theory}, CBMS Regional Conference Series
  in Mathematics, vol.~92, Published for the Conference Board of the
  Mathematical Sciences, Washington, DC, 1997. \MR{MR1421568 (97k:58183)}

\bibitem[Dod84]{dodziuk:84}
J.~Dodziuk, \emph{Difference equations, isoperimetric inequality and transience
  of certain random walks}, Trans. Amer. Math. Soc. \textbf{284} (1984), no.~2,
  787--794. \MR{MR743744 (85m:58185)}

\bibitem[EP05]{exner-post:05}
P.~Exner and O.~Post, \emph{Convergence of spectra of graph-like thin
  manifolds}, Journal of Geometry and Physics \textbf{54} (2005), 77--115.

\bibitem[EP07]{exner-post:07}
\bysame, \emph{{Convergence of resonances on thin branched quantum wave
  guides}}, to appear in Journal of Mathematical Physics (2007).

\bibitem[F07]{fulling.talk:07} S.~Fulling, \emph{Is there an
    interesting index theory for quantum graphs?}, talk at the Isaac
  Newton Institute (INI), Cambridge, (2007-02-27).
  
  \texttt{http://www.newton.cam.ac.uk/programmes/AGA/sem.html}

\bibitem[FKuW07]{fkw:pre07}
S.~Fulling, P.~Kuchment, and J.~H. Wilson, \emph{Index theorems for quantum
  graphs}, Preprint (2007).

\bibitem[FT04a]{friedman-tillich:pre04}
J. Friedman and J.-P. Tillich, \emph{Calculus on graphs}, Preprint
  \texttt{arXiv:cs.DM/0408028} (2004).

\bibitem[FT04b]{friedman-tillich:04}
\bysame, \emph{Wave equations for graphs and the edge-based {L}aplacian},
  Pacific J. Math. \textbf{216} (2004), no.~2, 229--266. \MR{MR2094545
  (2005k:05142)}

\bibitem[Gil95]{gilkey:95}
P.~B. Gilkey, \emph{{Invariance theory, the heat equation and the Atiyah-Singer
  index theorem}}, CRC Press, Boca Raton, 1995.

\bibitem[Har00]{harmer:00b}
M.~Harmer, \emph{Hermitian symplectic geometry and extension theory}, J. Phys.
  A \textbf{33} (2000), no.~50, 9193--9203. \MR{MR1804888 (2001m:58009)}

\bibitem[HP06]{hislop-post:pre06}
P.~Hislop and O.~Post, \emph{Exponential localization for radial random quantum
  trees}, Preprint (math-ph/0611022) (2006).

\bibitem[KPS07]{kps:pre07}
V.~Kostrykin, J.~Potthoff, and R.~Schrader, \emph{Heat kernels on metric graphs
  and a trace formula}, Preprint math-ph/0701009 (2007).

\bibitem[KS99]{kostrykin-schrader:99}
V.~Kostrykin and R.~Schrader, \emph{Kirchhoff's rule for quantum wires}, J.
  Phys. A \textbf{32} (1999), no.~4, 595--630.

\bibitem[KS03]{kostrykin-schrader:03}
\bysame, \emph{Quantum wires with magnetic fluxes},
  Comm. Math. Phys. \textbf{237} (2003), no.~1-2, 161--179, Dedicated to Rudolf
  Haag. \MR{MR2007178}

\bibitem[KS06]{kostrykin-schrader:06}
\bysame, \emph{Laplacians on metric graphs: eigenvalues,
  resolvents and semigroups}, Quantum graphs and their applications, Contemp.
  Math., vol. 415, Amer. Math. Soc., Providence, RI, 2006, pp.~201--225.
  \MR{MR2277618 (2007j:34041)}

\bibitem[Ku04]{kuchment:04}
P.~Kuchment, \emph{Quantum graphs: {I}. {S}ome basic structures}, Waves Random
  Media \textbf{14} (2004), S107--S128.

\bibitem[Ku05]{kuchment:05}
\bysame, \emph{Quantum graphs. {II}. {S}ome spectral properties of quantum and
  combinatorial graphs}, J. Phys. A \textbf{38} (2005), no.~22, 4887--4900.
  \MR{MR2148631}

\bibitem[KuZ01]{kuchment-zeng:01}
P.~Kuchment and H.~Zeng, \emph{Convergence of spectra of mesoscopic systems
  collapsing onto a graph}, J. Math. Anal. Appl. \textbf{258} (2001), no.~2,
  671--700.

\bibitem[KuZ03]{kuchment-zeng:03}
\bysame, \emph{Asymptotics of spectra of {N}eumann {L}aplacians in thin
  domains}, Advances in differential equations and mathematical physics
  (Birmingham, AL, 2002), Contemp. Math., vol. 327, Amer. Math. Soc.,
  Providence, RI, 2003, pp.~199--213. \MR{1 991 542}

\bibitem[L{\"u}c02]{lueck:02}
W.~L{\"u}ck, \emph{{$L\sp 2$}-invariants: theory and applications to geometry
  and {$K$}-theory}, Ergebnisse der Mathematik und ihrer Grenzgebiete.,
  vol.~44, Springer-Verlag, Berlin, 2002. \MR{MR1926649 (2003m:58033)}

\bibitem[MW89]{mohar-woess:89}
Bojan Mohar and Wolfgang Woess, \emph{A survey on spectra of infinite graphs},
  Bull. London Math. Soc. \textbf{21} (1989), no.~3, 209--234. \MR{MR986363
  (90d:05162)}

\bibitem[Nic87]{nicaise:87}
S.~Nicaise, \emph{Spectre des r\'eseaux topologiques finis}, Bull. Sci. Math.
  (2) \textbf{111} (1987), no.~4, 401--413. \MR{MR921561 (89a:58114)}

\bibitem[Ogu02]{ogurisu:02}
O.~Ogurisu, \emph{Supersymmetric analysis of the spectral theory on infinite
  graphs}, Seminars on infinite graphs and their spectrum at Lake Kawaguchi,
  January 2002, Contemporary Mathematics, 2002, pp.~57--75.

\bibitem[Pan06]{pankrashkin:06a}
K.~Pankrashkin, \emph{Spectra of {S}chr\"odinger operators on equilateral
  quantum graphs}, Lett. Math. Phys. \textbf{77} (2006), no.~2, 139--154.
  \MR{MR2251302 (2007f:81089)}

\bibitem[P05]{post:05}
O.~Post, \emph{Branched quantum wave guides with {D}irichlet boundary
  conditions: the decoupling case}, Journal of Physics A: Mathematical and
  General \textbf{38} (2005), no.~22, 4917--4931.

\bibitem[P06]{post:06}
\bysame, \emph{Spectral convergence of quasi-one-dimensional spaces}, Ann.
  Henri Poincar\'e \textbf{7} (2006), no.~5, 933--973. \MR{MR2254756}

\bibitem[Rot84]{roth:84}
J.-P. Roth, \emph{Le spectre du laplacien sur un graphe}, Th\'eorie du
  potentiel (Orsay, 1983), Lecture Notes in Math., vol. 1096, Springer, Berlin,
  1984, pp.~521--539.

\bibitem[RS01]{rubinstein-schatzman:01}
J.~Rubinstein and M.~Schatzman, \emph{Variational problems on multiply
  connected thin strips. {I}. {B}asic estimates and convergence of the
  {L}aplacian spectrum}, Arch. Ration. Mech. Anal. \textbf{160} (2001), no.~4,
  271--308. \MR{1 869 667}

\bibitem[Shi00]{shirai:00}
T. Shirai, \emph{The spectrum of infinite regular line graphs}, Trans.
  Amer. Math. Soc. \textbf{352} (2000), no.~1, 115--132. \MR{MR1665338
  (2001f:05092)}

\end{thebibliography}
\renewcommand{\MR}[1]{}

\end{document}